\documentclass[11pt,twoside]{article}
\usepackage{amsmath}
\usepackage{amsfonts}
\usepackage{amssymb}
\usepackage{graphicx}
\usepackage{epsfig}

\newcommand{\beq}[1]{\begin{equation}\label{#1}}
\newcommand{\eeq}{{\end{equation}}}
\newcommand{\beqn}[1]{\begin{eqnarray}\label{#1}}
\newcommand{\eeqn}{\end{eqnarray}}
\newcommand{\beaa}{\begin{eqnarray*}}
\newcommand{\eeaa}{\end{eqnarray*}}

\begin{document}

\newcommand{\ls}[1]
   {\dimen0=\fontdimen6\the\font \lineskip=#1\dimen0
\advance\lineskip.5\fontdimen5\the\font \advance\lineskip-\dimen0
\lineskiplimit=.9\lineskip \baselineskip=\lineskip
\advance\baselineskip\dimen0 \normallineskip\lineskip
\normallineskiplimit\lineskiplimit \normalbaselineskip\baselineskip
\ignorespaces }
\renewcommand{\theequation}{\thesection.\arabic{equation}}

\newtheorem{definition}[equation]{Definition}
\newtheorem{lemma}[equation]{Lemma}
\newtheorem{proposition}[equation]{Proposition}
\newtheorem{corollary}[equation]{Corollary}
\newtheorem{example}{Example}[section]
\newtheorem{conjecture}{Conjecture}[section]
\newtheorem{algorithm}{Algorithm}[section]
\newtheorem{theorem}[equation]{Theorem}
\newtheorem{exercise}[equation]{Exercise}
\newtheorem{assumption}[equation]{Assumption}

\def\le{\leq}
\def\ge{\geq}
\def\lt{<}
\def\gt{>}

\newcommand{\Exercise}{{\bf Exercise}}
\newcommand{\req}[1]{(\ref{#1})}
\newcommand{\beps}{{\scriptscriptstyle{{ \cal E}}}}
\newcommand{\lip}{\langle}
\newcommand{\rip}{\rangle}
\newcommand{\lf}{\lfloor}
\newcommand{\lc}{\lceil}
\newcommand{\rc}{\rceil}
\newcommand{\rf}{\rfloor}
\newcommand{\supp}{{\rm supp\,}}
\newcommand{\ARn}
 {\begin{array}[t]{c}
  \longrightarrow \\[-0.3cm]
 \scriptstyle {n\rightarrow \infty}
  \end{array}}

\newcommand{\uu}{\underline}
\newcommand{\oo}{\overline}
\newcommand{\dfr}{\displaystyle\frac}
\newcommand{\La}{\Lambda}
\newcommand{\la}{\lambda}
\newcommand{\eps}{\varepsilon}
\newcommand{\om}{\omega}
\newcommand{\Om}{\Omega}
\newcommand{\Inv}{{\text{\rm Inv}\,}}
\newcommand{\crit}{{\rm crit}}

\def\ot{{\overline \tau}}
\def\oX{{\overline X}}
\newcommand{\Si}{\Sigma}
\newcommand{\si}{\sigma}
\newcommand{\dist}{\rm dist}
\def\C{{\rm Const.}}
\def\leqsto{\leq_{\rm st}}


\newcommand{\EE}{{\mathbb E}}
\newcommand{\NN}{{\mathbb N}}
\newcommand{\PP}{{\mathbb P}}
\newcommand{\QQ}{{\mathbb Q}}
\newcommand{\reals}{{\mathbb R}}
\newcommand{\ZZ}{{\mathbb Z}}

\newcommand{\calA}{{\mathcal A}}
\newcommand{\calE}{{\mathcal E}}
\newcommand{\calB}{{\mathcal B}}
\newcommand{\calF}{{\mathcal F}}
\newcommand{\calG}{{\mathcal G}}
\newcommand{\calH}{{\mathcal H}}
\newcommand{\calI}{{\mathcal I}}
\newcommand{\calO}{{\mathcal O}}
\newcommand{\calQ}{{\mathcal Q}}
\newcommand{\calS}{{\mathcal S}}
\newcommand{\calV}{{\mathcal V}}
\newcommand{\calW}{{\mathcal W}}
\newcommand{\calY}{{\mathcal Y}}
\newcommand{\calX}{{\mathcal X}}

\newcommand{\bfcdot}{{\boldsymbol \cdot}}

\newcommand{\won}{{\boldsymbol 1}}
\newcommand{\hilbert}{\bigcirc\kern -0.8em
              {\rm\scriptstyle {H}\;}} 

\newcommand{\limn}{\lim_{n \rightarrow \infty}}
\newcommand{\limk}{\lim_{k \rightarrow \infty}}
\newcommand{\limi}{\lim_{i \rightarrow \infty}}
\newcommand{\liml}{\lim_{\ell \rightarrow \infty}}
\newcommand{\limv}{\lim_{v \rightarrow \infty}}
\newcommand{\limm}{\lim_{m \rightarrow \infty}}
\newcommand{\limd}{\lim_{\delta \rightarrow \infty}}
\newcommand{\limsupn}{\limsup_{n \rightarrow \infty}}
\newcommand{\liminfn}{\liminf_{n \rightarrow \infty}}

\newcommand{\proof}{\noindent \textbf{Proof:\ }}
\newcommand{\remark}{\noindent \textbf{Remark:\ }}
\newcommand{\remarks}{{\bf Remarks:}}

\newcommand{\eff}{{\operatorname{eff}}}

\newcommand{\dsum}{\displaystyle\sum}
\newcommand{\dprod}{\displaystyle\prod}

\newcommand{\ffrac}[2]
  {\left( \frac{#1}{#2} \right)}

\newcommand{\one}{\frac{1}{n}\:}
\newcommand{\half}{\frac{1}{2}\:}

\def\squarebox#1{\hbox to #1{\hfill\vbox to #1{\vfill}}}
\newcommand{\qed}{\hspace*{\fill}
           \vbox{\hrule\hbox{\vrule\squarebox{.667em}\vrule}\hrule}\smallskip}


\title{\textbf{A law of large numbers for random walks in random 
mixing environments}}

\author{Francis~Comets%
 \thanks{
Universit{\'e} Paris 7 --- Denis Diderot,
Math{\'e}matiques, case 7012,
2 place Jussieu,
752151 Paris Cedex 05,
France.}
\and Ofer~Zeitouni%
 \thanks{Departments of Electrical Engineering and of Mathematics, Technion,
Haifa 32000, Israel. Partially supported by the Technion V.P.R. fund for
the promotion of research.
Part of this work was done while the author visited
Univ. Paris 7.}}

\date{\large{\textit{ 29  April 2002}}}
\maketitle
\pagestyle{myheadings}
\markboth{{\sc LLN for RW in mixing RE}}{{\sc Comets and Zeitouni}}

\begin{abstract}
We prove a law of large numbers for a class of multidimensional
random walks
 in random environments where the environment satisfies appropriate
 mixing conditions, which hold when the environment is a
weak mixing  field
in the sense of Dobrushin and Shlosman.
Our result  holds if the mixing 
rate balances moments of some random times
depending on the path. It applies in the non-nestling case, but 
we also provide examples of nestling walks that satisfy
our assumptions. The derivation is based on an adaptation,
using coupling, of the regeneration
argument of Sznitman-Zerner.

\end{abstract}

\noindent
{\small {\bf Key Words:} Random walk
 in random environment, law of large numbers, Kalikow's condition,
nestling walk, mixing.
\\
{\bf AMS (1991) subject classifications:} 60K40, 82D30.\\
{\bf Short title:} Random walk in mixing environment }
        
\section{Introduction and statement of results}
Let $S$ denote the 2d-dimensional simplex, and
set $\Om=S^{\ZZ^d}$. We consider $\Om$ as 
an ``environment'' for the random walk that we define below in (\ref{1}).
We denote by $\om(z,\cdot)=\{\om(z,z+e)\}_{e\in \ZZ^d, |e|=1}$ the 
coordinate of $\om\in\Om$ corresponding to $z\in \ZZ^d$.

Conditional on a realization $\om\in\Om$, 
we define the Markov Chain $\{X_n\}=\{X_n; n \geq 0\}$ with
state space $\ZZ^d$ started at $z\in \ZZ^d$
as the process satisfying $X_0\equiv z$ and
\begin{equation}
\label{1}
P_\om^z (X_{n+1}  = x+ e | X_n = x) = \om (x,x+e)\,,\quad e\in  \ZZ^d, 
|e|=1\,.
\end{equation}
The law of the random walk in random environment (RWRE)
$\{X_n\}$  under this transition kernel, denoted 
$P_\om^z(\cdot)$, depends on the environment $\om\in\Om$ and is called the
\textit{quenched} law of $\{X_n\}$. 

Let $P$ be a probability measure on $\Om$, stationary and ergodic
with respect to the shifts in  $\ZZ^d$.
With a slight abuse of notations, we
write
$\PP^z= P\otimes P_\omega^z$
for both the joint law on $\Om \times (\ZZ^d)^\NN$ of
$\{X_n\}_n$ and $\om$, and for its marginal on
$(\ZZ^d)^\NN$; in the latter case, we refer to it as the \textit{annealed }
law of the  process $\{X_n\}$. We will denote
by $\EE^z=E_{\PP^z}, E_\om^z=E_{P_\om^z}$ the expectations
corresponding to $\PP^z, P_\om^z$, respectively. Considering the
 annealed law rather than the quenched law, one takes  advantage of 
the smoothing from the $\om$-average, 
but the Markov property is lost.

The RWRE with $d=1$ is by now well studied, see \cite{stflour}
for a recent review. The multidimensional case is much less
understood. A crucial simplification in the case $d=1$ is that a nearest neighbor 
random walk tending to $+ \infty$ has to visit all positive sites; then one 
can use ergodic theorem to smooth  the environment out. 
In constrast, it is not clear how to take advantage of ergodicity of the medium
in dimension $d >1$.
When $P$ is a product measures,
laws of large numbers  and central limit theorems 
for $\{X_n\}$ were derived in an 
impressive
sequence of papers \cite{sznitmanzerner}, \cite{Sz1}, \cite{Sz2},
focusing on the ballistic regime.
Our goal in this paper is to present a technique, based on an appropriate
coupling,  for
extending some of the results of \cite{sznitmanzerner},
\cite{Sz1} 
to the case where $P$ is not a product measure.
 
A motivation for this question 
relates to an example in  \cite[Proposition 2]{zernermerkl}, of a 
RWRE in dimension 2
with the non-standard asymptotics: $\PP^o( \lim_n X_n/n = w)=1/2$, 
$\PP^o( \lim_n X_n/n = -w)=1/2$ for some 
non zero vector $w$. There, the environment is ergodic (but not mixing),
and the RWRE
is strictly elliptic (but not uniformly). In view of this example, it seems
important to clarify the specific role of the various assumptions
used to get a standard law of large numbers, e.g. independent,
identically distributed (i.i.d.) environment,
uniform ellipticity and drift-condition in
\cite{sznitmanzerner}.

We work in the context of {\it ballistic} walks, i.e. walks
$X_n$ which 
tend to infinity in some direction $\ell\in \reals^d\setminus \{0\}$,
with a non-vanishing speed.
Conditions for the first statement to occur have been explored by
Kalikow \cite{kalikow} two decades ago, see below Assumption 
$({\cal A}4)$.
With an i.i.d. environment, Sznitman and Zerner \cite{sznitmanzerner}
introduced a sequence of {\it regeneration} 
times
and showed, roughly, that the environments traversed by the walk between
regeneration times, together with the path of the walk, form a sequence of 
i.i.d. random vectors under the {\it annealed} law $\PP^o$.
This allowed them to 
derive a law of large numbers under Kalikow's condition
by studying the tail properties of these regeneration times.

In \cite{stflour}, 
a {\it coupling technique} was introduced 
that immediately allows one to adapt the construction of regeneration
times to the
case of 
measures $P$  which are $L$-dependent, that is such that coordinates 
of the environment at distance larger than $L$, some fixed deterministic
$L$,  are independent. This coupling covers in particular the setup in
\cite{shen}, that deals with a particular 1-dependent environment
(note however that we do not attempt to
recover here all the results of \cite{shen}).
Intuitively, the coupling idea is that, due to the
uniform ellipticity property,
the walk has positive probability for travelling the  $L$ next steps
whithout looking at the environment. One then is reduced to the study of 
tails of these regeneration times. 

Our purpose in this work is to further modify 
the construction of regeneration times
and allow for more general type of mixing conditions on $P$.
A complication arises from the destruction of the
renewal structure due to the dependence in the environment: in fact,
the environments between regeneration times need not even define a stationary
sequence any more. Our approach is based on suitably approximating
this sequence by an i.i.d. sequence, namely the so-called splitting 
representation \cite{thorisson}.

Another  approach to non product measures $P$, but with a rather mild
dependence structure, has been proposed in 
\cite{komorowski}, \cite{komorowski2}. A comparison between their 
results and ours is presented at
the end of the article. Our results here cover natural examples of
environment distributions as Gibbs measures in a mixing regime.

We now turn to a description of our results.
We deal with environments subject to various mixing 
conditions, and it is appropriate to consider closed positive {\it cones}.
For  $\ell \in\reals^d \setminus \{0\}$, $x \in \reals^d$ and
$\zeta \in (0,1)$, 
define the cone of vertex $x$, direction $\ell$ and angle
$\cos^{-1} ( \zeta)$,
\beq{defcone}
C(x, \ell, \zeta)=\{y \in \reals^d \;; \;
(y-x) \cdot \ell \geq \ \zeta |y-x| |\ell| \}\;.
\end{equation}
(All through the paper, $|\cdot|$ denotes  the  euclidean
norm on $\reals^d$, and $|\ell|_1=\sum_i |\ell_i|$ the $\ell_1$-norm.)
Note that for $\zeta=0$, this is just a usual half-space.

In the sequel, we fix an 
$\ell\in \reals^d\setminus \{0\}$ such that $\ell$
has {\it integer coordinates}. With ${\rm sgn}(0)=0$, let
\beq{calE}
\calE_{\bar \eps}=\{{\rm sgn}(\ell_i)e_i\}_{i=1}^d
\setminus \{0\}\;.
\end{equation}

Throughout, we make the following two assumptions
on the environment:
\begin{assumption}
\label{ass-RWRE}
\mbox{}
\begin{itemize}
\item[(${\cal A} 1$)]
$P$ is stationary and ergodic, and satisfies the following mixing 
condition on $\ell$-cones: 
for all positive $\zeta$ small enough
there exists a function
$\phi (r) \underset{r\to\infty}{\to} 0$ such that any two
events $A, B$ with $P(A) > 0$, $A \in \si\{ \om_z; z \cdot \ell
\leq 0\}$ and $B   \in \si\{ \om_z; z \in 
C(r \ell, \ell, \zeta)\}$ 
it holds that
$$
\left| \frac{P(A \cap B)}{P(A)} - P(B) \right| \le \phi (r |\ell|)\,.
$$
\item[(${\cal A} 2$)]
$P$ is elliptic and { uniformly elliptic with respect to} $\ell$:
$P(
\om (0,e)>0 ; |e|=1)=1$, and 
there exists a $\kappa>0$ such that
$$
P(\min_{e\in \calE_{\bar\eps}}
\om (0,e)\geq \kappa)=1
\,.
$$
\end{itemize}
\end{assumption}
Note that $({\cal A} 1)$ is equivalent to
\beq{fmix}
|E (fg) -Ef Eg| \leq \phi (r |\ell|) \|f\|_1 \|g\|_{\infty}
\end{equation}
for all bounded functions, $f$ being $\si\{ \om_z; z \cdot \ell
\leq 0\}$-measurable, and $g$ being $\si \{ \om_z; z \in
C(r \ell, \ell, \zeta)\}$-measurable.
Properties of the type $({\cal A} 1)$ are generically called $\phi$-mixing or
uniform mixing \cite[Section 1.1]{doukhan}.
When (${\cal A} 1$) holds for $\zeta=0$, we will say, in this paper, 
that the field $P$ is $\phi$-mixing. But this condition is too
restrictive,
and we give in Section \ref{sec:mix}, 
examples of environments satisfying
the mixing property $({\cal A} 1)$, but not $\phi$-mixing (with
$\zeta=0$). 
In most applications, one can find 
such a $\phi$ not depending on 
 $\ell \in\reals^d \setminus \{0\}$, for which  the condition holds for all  $\ell$.

Next, we turn to conditions on the environment ensuring the ballistic nature
of the walk. In order to do so, we introduce (a natural extension of) 
Kalikow's Markov chain \cite{kalikow} as follows.
Let $U$ be a finite, connected subset of $\ZZ^d$, with $0\in U$, let
$$
\calF_{U^c}=\sigma\{\om_z: z\not\in U\}\;,
$$
and define on $U \cup\partial U$ an auxiliary Markov chain with transition
probabilities
\begin{equation}
\label{kal1}
\hat P_U (x,x+e) = \begin{cases}
\frac{\displaystyle \EE^o\left[\sum_{n=0}^{T_{U^c}}\;\won_{\{X_n=x\}}
\om(x,x+e)|\calF_{U^c}\right]}
{\displaystyle \EE^o \left[\sum_{n=0}^{T_{U^c}} \:\won_{\{X_n=x\}}
|\calF_{U^c}\right]}, & x\in U,  |e|=1\\
1& x \in \partial U, e=0
\end{cases}
\end{equation}
where $T_{U^c} = \min \{n \ge 0: X_n \in \partial U\}$ (note that the
expectations in (\ref{kal1}) are finite due the Markov property
and $\ell$-ellipticity).
The transition kernel $\hat P_U$ weights the transitions
$x\mapsto x+e$ according to the occupation time of the vertex $x$ before
exiting $U$.  
Define the Kalikow drift as
$\hat d_U(x)=\sum_{|e|=1}  e
\hat P_U (x,x+e)$, with the RWRE's drift at
$x$ defined by $d(x,\om)=\sum_{|e|=1}  e \om(x,x+e)$. 
Note that, unlike in the i.i.d. case,
 the Kalikow drift, as well as Kalikow's chain itself,  here is random
because it depends on the environment outside of $U$.
This new  Markov chain is useful because of the 
following property \cite{kalikow}, which remains valid in our non-i.i.d. setup.
Since $U$ is finite and the walk is uniformly elliptic in the direction
$\ell$, under both
$\hat P_{ U}$ and 
$\PP^o( \cdot |\calF_{U^c})$, the exit time $T_{U^c}=\inf\{n \geq 0;
X_n \in U^c\}$ is finite, 
and 
\begin{equation}
\label{sortie}
X_{T_{U^c}} {\rm \ has\ the\ same\ law\ under\ }  \hat P_{ U} {\rm \ and\ } 
\PP^o( \cdot |\calF_{U^c})\;.
\end{equation}
In this paper, we will consider one of
following {\it drift} conditions, which ensure 
a ballistic behavior for the walk:
\begin{assumption}
\label{ass-RWRE1}
\mbox{}
\begin{itemize}
\item[(${\cal A} 3$)] Kalikow's condition: There exists a $\delta(\ell)>0$
deterministic such that
$$
\inf_{U, x \in U} \hat d_U(x) \cdot \ell \geq \delta(\ell), \, P-a.s..
$$
(The infimum is taken over all connected finite subsets
of $\ZZ^d$ containing 0.)
\item[(${\cal A} 4$)] Non-nestling: There exists 
  a $\delta(\ell)>0$ such that
$$
 d(x, \om) \cdot \ell \geq \delta(\ell)\,, P-a.s.
$$
\end{itemize}
\end{assumption}

We will always assume $({\cal A} 1)$ and (${\cal A} 2$), 
and (except in the beginning of Section \ref{sec:LLN}), also
one of (${\cal A} 3$) or the stronger (${\cal A} 4$).
Clearly,  (${\cal A} 4$) implies that  (${\cal A} 3$) holds with the
same $\ell$ and $\delta(\ell)$.
An inspection of the proof in \cite{kalikow} reveals that, 
under  (${\cal A} 3$), the conclusion
$\PP^o (\liminf_{n\to\infty} X_n\cdot\ell=\infty
)=1$ remains valid in our non-i.i.d. setup
(see e.g. the exposition in \cite{stflour}).  Note that the requirement 
from $\ell$ to possess integer coordinates is not a restrictive one
for given an $\ell\in \reals^d\setminus\{0\}$ satisfying
either (${\cal A} 3$) or (${\cal A} 4$), one 
finds by continuity an $\ell$
with integer coordinates satisfying the same. 
We make this restriction for the 
convenience of defining the path $\bar \eps$ 
in (\ref{stayincone}) below.

The statement of our fundamental result, Theorem \ref{theo-LLN}, 
involves 
certain {\it modified regeneration times}, which are introduced in
Section \ref{sec:times} below. The basic condition in this theorem
is a trade-off between moments for these ``regeneration times''
and rate of mixing for the environment.
A corollary can be readily stated 
here: it deals with 
the non-nestling case, and does not require any assumption on the
mixing rate other than $\phi \to 0$.
\begin{corollary}
\label{cor-nonnest}
Assume $({\cal A}1, {\cal A}2)$ and $({\cal A}4)$ hold for some 
$\ell \in \ZZ^d\setminus \{0\}$.
 Then there exists a deterministic
$v\neq 0$  such that 
$$\lim_{n\to\infty}
\frac{ X_n}{n}=v\,,\quad \PP^o-a.s..$$
\end{corollary}
Our fundamental Theorem \ref{theo-LLN} does not restrict 
attention to non-nestling cases,
but applies also when only Kalikow's condition  holds.
A class of nestling examples based on Ising-like environments, 
satisfying the conditions of Theorem
\ref{theo-LLN} and hence the law of large numbers, 
is provided in Theorem \ref{theo-exgibbs} below.
Dealing with nestling walks is much more delicate than with non-nestling
ones due to the existence of so-called {\it traps},
i.e., finite but large regions where the environment is atypical,
confining the walk or creating an abnormal drift. In this case, a renormalization
procedure (``coarse graining'') is needed to control the size of the traps 
and their effect.

The structure of the article is as follows: in Section \ref{sec:times},
we introduce the coupling representation (\ref{numeps})
which allows us to deal with non-independent
medium, and some {\it approximate regeneration times}
$\tau_i^{(L)}$ parametrized by a parameter $L$: they are defined by
(\ref{D'}) and (\ref{num8}), and they lead to an ``approximate renewal''
result, our Lemma \ref{lem-"renewal"} below. In
Section \ref{sec:LLN}, 
we prove  the law of large numbers under $({\cal A}1,2,3)$ and a
suitable integrability condition $({\cal A} 5)$. We then give the proof 
of Corollary \ref{cor-nonnest}, by showing that Condition $({\cal A} 5)$
is trivially satisfied under the non-nestling assumption $({\cal
A}4)$. In Section \ref{sec:mix},
we make a short digression
to show that the mixing assumption $({\cal A}1)$ is satisfied in many 
cases of interest.  In Section \ref{sec-A5} we show how 
condition $({\cal A} 5)$ can be checked in the nestling case,  
and we  construct  a class of nestling, mixing environments
satisfying our conditions for the law
 of large numbers, see 
Theorems \ref{theo-A5} 
and \ref{theo-exgibbs} for precise statements.
Finally, Section \ref{sec:conclusions} is devoted to
concluding remarks and extensions. 
 
\section{Some Random Times}  \label{sec:times}
\setcounter{equation}{0}

To implement our coupling technique, 
we begin, following \cite{stflour}, by constructing an extension of 
the probability space, depending on the vector
$\ell$ with integer coordinates: recall that 
the RWRE was defined by means of the law
$\PP^o = P\otimes P_\om^o$ on the canonical space 
$(\Om \times (\ZZ^d)^{\NN}, \calF \times \calG)$.
Set $W = \{0\}\cup \calE_{\bar \eps}$ (recall
the notation (\ref{calE})),
and let $\calW$ be the cylinder
$\sigma$-algebra on $W^\NN$.
We now define the measure
$$
\oo\PP^o = P \otimes Q \otimes \oo{P}^o_{\om,{\beps}}
$$
on
$$
\Bigl( \Om \times W^\NN \times (\ZZ^d)^\NN, \quad \calF \times \calW \times
\calG\Bigr)
$$
in the following way: 
$Q$ is a product measure, such that with $\beps
=(\eps_1,\eps_2,\ldots)$
denoting an element of $W^\NN$,
$Q (\eps_1 =  e) = \kappa$, for $e\in \calE_{\bar \eps}$,
while
$Q(\eps_1=0) = 1-\kappa |\calE_{\bar\eps}|$.
For each fixed $\om,\beps$, $\oo{P}_{\om,\beps}^o$ is the 
law of the Markov chain $\{X_n\}$ with state space $\ZZ^d$, such that
$X_0=0$ and, for every $z, e\in \ZZ^d$, $|e|=1$,
\begin{equation} 
\label{numeps}
\oo{P}_{\om,\beps}^o (X_{n\!+\!1} = z\!+\!e \;|\; X_n \!=\! z) =
\won_{\{\eps_{n\!+\!1} =e\}} +
\frac{\won_{\{\eps_{n\!+\!1}=0\}}}{1-\kappa |\calE_{\bar\eps}|}
[\om (z,z\!+\!e)-\kappa  {\bf 1}_{\{e\in \calE_{\bar\eps}\}}]\,.
\end{equation} 
Clearly, the law of $\{X_n\}$ under 
$Q \otimes \oo{P}_{\om,\beps}^o$ coincides
with its law under $P_\om^o$, while its law under
$\oo{\PP}^o$ coincides with its law under $\PP^o$.

We fix now a particular sequence of
$\beps$ in $\calE_{\bar \eps}$ of length
$|\ell|_1$  with sum  equal to $\ell$:
for definiteness, we take $\bar \eps = (\bar \eps_1, \ldots,
\bar \eps_{|\ell|_1})$ with 
\beaa
\bar \eps_1= \bar \eps_2= \ldots \bar \eps_{|\ell_1|}= {\rm sgn}
(\ell_1) e_1\;, \;
\bar \eps_{\ell_1+1}= \bar \eps_{\ell_1+2}= \ldots
 \bar \eps_{|\ell_1|+|\ell_2|}= {\rm sgn} (\ell_2) e_2\;,  \\
\ldots \;\; \bar \eps_{|\ell|_1-|\ell_d|+1}=  \ldots
 \bar \eps_{|\ell|_1}= {\rm sgn} (\ell_d) e_d\;.
\eeaa
We fix, from now on through the whole paper, $\zeta>0$ small enough such that
\begin{equation}
\label{stayincone}
\bar \eps_1, \bar \eps_1+\bar \eps_2, \ldots, \bar \eps_1+
\ldots
\bar \eps_{|\ell|_1}=\ell \in C(0, \ell, \zeta)\;\;,
\end{equation}
and such that $({\cal A}1)$ above is satisfied. Without mentioning it
explicitely in the sequel, we always rotate the 
axes such that $\ell_1\neq 0$.

For $L \in |\ell|_1 \NN^*$ we will denote by $\bar \eps^{(L)}$ the vector
$$
\bar \eps^{(L)}=(\bar \eps, \bar \eps, \ldots, \bar \eps)
$$
of dimension $L$. In particular, $\bar \eps=\bar \eps^{(|\ell|_1)}$,
and for $\beps$ with $(\eps_{n+1}, \ldots \eps_{n+L})=\bar \eps^{(L)}$,
$$
\oo{P}_{\om,\beps}^o (X_{n+L} = x + \frac{L}{|\ell|_1} \ell \; \vert
\; X_{n}=x) =1\;\;,
$$
and the path $X_n, X_{n+1}, \ldots X_{n+L}$ remains in the cone
$C(x, \ell, \zeta)$.

Define 
\begin{equation} 
\label{D'}
D'=\inf\{n\ge 0:\, X_n \notin C(X_0, \ell , \zeta)\}\;\;.
\end{equation}
(For $\zeta=0$, this is $D$ from \cite{sznitmanzerner}.) Let us
state a few direct consequences of  (${\cal A} 3$) and  (${\cal A} 4$).

\begin{lemma}
\label{D'>0} Assume (${\cal A} 3$). Let $f(y)=y \cdot \ell - \zeta 
|y| |\ell|$
with $\zeta \leq \delta(\ell)/(3|\ell|)$.

1) There  exist a $\lambda_0=\lambda_0(\delta(\ell))>0$ such that for
all $\lambda \in (0, \lambda_0]$ and 
all connected, finite subset $U$ of $\ZZ^d$ containing
0, 
$$ 
M_n^{\lambda}=\exp \{-3 \la f(X_n) +  \lambda \delta(\ell)
 (n\wedge T_{U^c})\}
$$ 
is a supermartingale for 
 Kalikow's Markov chain $\hat P_{ U}$ from (\ref{kal1})
($T_{U^c}$ is the exit time of
 $U$ for $ X_n$).

2) For $m>|\ell |$, consider
the truncated cone $ V_m=C(0, \ell, \zeta) \bigcap \{ y \in \reals^d
; y \cdot \ell \leq m \}$.
We have $\hat P_{ V_m}(X_{T_{V_m^c}} \cdot \ell > m) \geq 2\eta$
with some constant $\eta>0$ depending on $\delta(\ell), \zeta$ 
but not  on $m, \om$.

3) If also 
(${\cal A} 4$) holds (and $\ell$ is general),  choosing
$\kappa$ small enough such that $\delta(\ell)>2 \kappa$, it holds
that
$M_n^{\lambda}$ is a supermartingale
under the quenched measure $\oo{P}_{\om,\beps}^o$
for all  $\om,\beps, \lambda 
\in (0, \lambda_0]$,
and $\oo{ P}_{\om, \beps}^o (X_{T_{V_m^c}} \cdot \ell > m) \geq
2\eta$. Moreover, with $W_r=\{ x \in \ZZ^d; x \cdot \ell < r\}, r>0$,
\begin{equation} 
\label{expexit-q}
\oo{ E}_{\om, \beps}^o ( \exp\{ \la \delta(\ell) T_{W_r^c}\} )\leq 
\exp\{3 \la r\}\;.
\end{equation}

4)  Assume (${\cal A} 3$) holds with $\ell = e_1$, and $\delta=\delta (e_1)$. 
Then, there exists $\la_1= \la_1(
\delta )>0$ with $\la_1 \to +\infty$ as $\delta \to 1^-$ 
such that  $\exp \{ -3 \la f(X_n)\}$ is a supermartingale for 
 Kalikow's Markov chain ($\la \in [0, \la_1)$). In particular, 
 $\inf_{m, \om} \hat P_{ V_m}(X_{T_{V_m^c}} \cdot \ell > m) \to 1$ 
when $\delta \to 1^-$.
\end{lemma}

We stress that the above constants do not depend on $\om$ outside $U$,
in contrast to Kalikow's Markov chain $\hat P_{ U}$ itself. Recall
that, due to (\ref{sortie}), estimates on the exit
distribution for Kalikow's Markov chain yields the similar estimate
for the RWRE, but on the other hand,  exit time distribution for Kalikow's Markov
chain and RWRE may be quite different. 

\begin{proof}
1) 
Since the chain has unit jumps,  (${\cal A} 3$) implies for
$y \in U$
$$
\hat E_{U}[(f(X_{n+1})-f(X_n)) | X_n=y] \geq (2/3) \delta(\ell)
$$
and $f(X_{n+1})-f(X_n)$ is uniformly bounded. 
Choosing $\la_0$ such that 
$\max\{u^{-2}(e^u-1-u); 0< |u|\leq \la_0 (3|\ell|+2)\} \leq 
\delta(\ell)/( \la_0[3|\ell|+2]^2)$, we have for $\la \in [0, \la_0]$
that
$\hat E_{ U} ( \exp\{- 3 \la (f(X_{n+1})-f(X_n)) 
 +  
\lambda \delta(\ell)\}| X_n=y) \leq 1 $ uniformly in  $y \in U$.
Note that we can choose $\la_0$ increasing in $\delta( \ell)$.

2) Applying the stopping theorem for the exit time $T_{V_m^c}$
to the above supermartingale,
we get for $y \in V_m$,
$$
\exp\{ -3\la f(y)\} \geq \hat E_{ V_m} ( \exp\{-3\la
f(X_{T_{V_m^c}})+ \lambda \delta(\ell)
 T_{V_m^c}   \} | X_0=y) \geq  \hat P_{ V_m}
(X_{T_{V_m^c}} \cdot \ell \leq m  | X_0=y)\;,
$$
where the second inequality is due to
$f<0$ on the boundary of the cone $C(0, \ell, \zeta)$ and 
$\la \delta(\ell) >0$. 
With $y=\bar \eps_1 \in  V_m$, we have, for the chain starting from 0,
\begin{align}
\label{newf}
\hat P_{ V_m}(X_{T_{V_m^c}} \cdot \ell > m) \geq 
\hat P_{ V_m}(X_1= \bar \eps_1 ) \hat P_{ V_m}
(X_{T_{V_m^c}} \cdot \ell \geq m  | X_1=\bar \eps_1) 
\geq \kappa [1-e^{-3 \la f( \bar \eps_1 )}]
\end{align}
which is positive since $f$ is positive in the interior of the cone.

3) It is straightforward to check that the above computations apply
to $\oo{ P}_{\om, \beps}^o$ under the assumption
 (${\cal A} 4$) (the assumption $\delta(\ell)>2\kappa$ is used
to ensure that modified environment appearing in the
definition (\ref{numeps})
of $\oo P_{\om,\beps}^o$ also is uniformly elliptic
in the direction $\ell$). 
In addition,
to prove (\ref{expexit-q}), we apply the stopping theorem to the
$\oo{ P}_{\om, \beps}^o$-supermartingale $M_n^{\lambda}$ at the exit
time of the domain $W_r$ intersected with large, finite boxes, and we
get
$$
\exp\{-3 \la r\} \oo{ E}_{\om, \beps}^o ( \exp\{ \la \delta(\ell) T_{W_r^c}\} )\leq 1
\;.
$$

4) With $\delta=\delta(e_1)$,
set ${\cal A}=\{(\alpha,\beta)\subset [0,1]^2: 
\alpha+\beta\leq 1, \alpha-\beta\geq \delta\}$.
Then, for 
 $y \in V_m$, and $\la>0$,
it follows from $({\cal A}3)$ that
\begin{eqnarray*}
\hat E_{U}[ \exp\{- 3 \la ((X_{n+1}-X_n)\cdot e_1 + 3 \la \zeta\} |  X_n=y] 
&\leq &
\exp\{3\la \zeta \} 
\sup_{(\alpha,\beta)\in {\cal A}}
\left[\alpha e^{-3\la}+\beta e^{3\la}+ (1-\alpha-\beta)\right]
\\
&=&
\exp\{3\la \zeta \} 
\left[ \cosh(3\la)-\delta \sinh(3\la)\right]=: A(\delta,\la,\zeta)\,.
\end{eqnarray*}
We see that $\la_1(\delta) := \sup\{\la:  A(\delta, \la, \zeta) <1 , 
\; \forall \zeta \leq \delta/3 \} 
\to \infty$ as $ \delta \to 1$. In particular, the right-hand
 side of (\ref{newf})
can be arbitrary close to 1 as  $ \delta \to 1$.
\qed
\end{proof}

Assumption $({\cal A} 3)$ implies $\PP^o (D'=\infty) > 0$. Indeed, consider
the truncated cone $
V_m=C(0, \ell, \zeta) \bigcap \{ y \in \reals^d ; y \cdot \ell \leq m \}
$, and Kalikow's Markov chain $\hat P_{ V_m}$. From (\ref{sortie}), 
the exit distribution for $X$ out of $V_m$ is the same 
under both 
$\hat P_{ V_m}$ and $ \PP^o(\cdot | \calF_{V_m^c})$.  From  part 2)
 of Lemma 
\ref{D'>0}, it follows
$$
\PP^o (D'=\infty |\; \om_x, x \cdot \ell \leq 0)=
\lim_{m \to \infty}  \PP^o (X_{T_{V_m^c}} \cdot \ell > m 
|\; \om_x, x \cdot \ell \leq 0) \geq 2 \eta\,, \quad P-a.s..
$$
By integration, we get
\begin{equation} 
\label{boundD'}  
\PP^o (D'=\infty\; |\; \om_x, x \cdot \ell \leq -r) \geq 
2 \eta, \;\;\;\; P-a.s. 
\end{equation}

Recall that $({\cal A} 3)$ 
implies that $\PP^o( X_n\cdot\ell\to_{n\to\infty}+\infty)=1$.
Set 
$$
\calG_n = \sigma ((\eps_i, X_i),\,i\le n)\;,
$$
 fix $L \in |\ell|_1
\NN$ 
and, setting $\oo{S}_0 = 0, 
$ define,
using $\theta_n$ to denote  time shift 
and $\bar \theta_x$ to denote space shift,
\begin{align} \nonumber
&
\oo{S}_1  =  \inf  \Bigr\{ n \geq  L: X_{n-L} \cdot \ell > 
\max \{X_m \cdot \ell:
\, m \!<\! n\!-\!L\},\,
(\eps_{n\!-\!1}, \dots, \eps_{n\!-\!L})= \bar \eps^{(L)} \Bigr\} \le 
\infty \;,\\
\label{num8}
&
 \oo{R}_1 = D' \circ \theta_{\bar{S}_1}
+ \oo{S}_1 \le \infty.
\end{align} 
Define further, by induction for $k\ge 1$,
\begin{align*} 
&
\oo{S}_{k+1}  =
\inf  \Bigr\{ n \ge R_k: X_{n-L} \cdot \ell > \max \{X_m \cdot \ell:
\, m\! <\! n\!-\!L\},\,
(\eps_{n\!-\!1}, \dots, \eps_{n\!-\!L})= \bar \eps^{(L)} \Bigr\} \le \infty
\;,\\
&
\oo{R}_{k+1} = 
D' \circ \theta_{\bar{S}_{k+1}}
+ \oo{S}_{k+1} \le \infty\;,
\end{align*} 

Clearly, these are $\calG_n$ stopping times (depending on $L$), and 
$$
0 = \oo{S}_0 \le \oo{S}_1 \le \oo{R}_1 \le \oo{S}_2 \le \cdots \le \infty
$$
and the inequalities are strict if the left member is finite.  On the
set $A_{\ell}:=\{X_n\cdot\ell \to_{n\to\infty} \infty\},$ it is straightforward
to check (using the product structure of $Q$), 
that the time $\oo{S}_1$ is $\oo \PP^o$-a.s. finite, as is
 $\oo{S}_{k+1}$ on the set $A_{\ell} \bigcap \{ \oo{R}_k < \infty\}$.
Define:
\begin{align*}
K & = \inf \{k \ge 1: \oo{S}_k < \infty, \oo{R}_k = \infty\} \le \infty,\\
\tau_1^{(L)} & = \oo{S}_K \le \infty\,.
\end{align*}
This random time $\tau_1^{(L)}$ is the first time $n$ when the walk performs as follow: at
time $n-L$ it has reached a record value in the direction $+\ell$,
then 
it travels using the $\beps$-sequence only up to time $n$, and from
time $n$ on, it doesn't exit the positive cone $C(X_n, \ell, \zeta)$
with vertex $X_n$. 
In particular, $\tau_1^{(L)}$ is not a stopping time, and 
we emphasize its dependence on $L$. As in \cite{stflour}, the advantage
in working with $\tau_i^{(L)}$ (as opposed to the more standard
$\tau_i^{(0)}$) is that the $\{\beps\}$ sequence creates a 
spacing where no
information on the environment is gathered by the RWRE.

\begin{lemma}
\label{finitetau} 
Assume $({\cal A}1,2,3)$, and  $\zeta \leq \delta(\ell)/(3|\ell|)$. 
Then, there exists a 
$L_0$ such that for $L \geq L_0$,
  $\tau_1^{(L)}$ is finite $\oo{\PP}^o$-a.s..
\end{lemma}

\begin{proof}
This amounts to proving $K < \infty$. Toward this end, write
\begin{align}
\oo\PP^o (\oo{R}_{k+1} & < \infty) = \oo\PP^o (\oo{R}_k < \infty, D' \circ 
\theta_{{\oo{S}_{k+1}}} < \infty) \nonumber \\
& = \sum_{z\in\ZZ^d}\sum_{n\in\NN} \oo\PP^o ( \oo{R}_k< \infty, D'\circ \theta_n < \infty,
X_{\oo{S}_{k+1}} = z, \oo{S}_{k+1}=n) 
\label{taufinite1}\\
& = \sum_{z\in\ZZ^d} \sum_{n\in\NN} E_{P\otimes Q}
\Bigl(\oo{P}_{\om,\beps}^o ( \oo{R}_k< \infty,
X_{\oo{S}_{k+1}} = z, \oo{S}_{k+1}=n, D'\circ \theta_n < \infty) 
 \Bigr)
 \nonumber \\
& = \sum_{z\in\ZZ^d} \sum_{n\in\NN} E_{P\otimes Q}
\Bigl(\oo{P}_{\om,\beps}^o (\oo{R}_k < \infty, X_{\oo{S}_{k+1}} = z, 
\oo{S}_{k+1} = n) 
\cdot \oo{P}^o_{{\bar \theta}_z{\om,\theta_n \beps}} (D' < \infty) \Bigr)
 \nonumber \,
\end{align}
using the strong 
Markov property under $\oo{P}_{\om,\beps}^o$. The point here is 
that $\oo{P}_{{\bar \theta}_z \om, \theta_n \beps}^o (D' < \infty)$ is
measurable on 
$\sigma(\om_x:\; x \in C(z, \ell, \zeta)) \otimes 
\sigma(\eps_i, i\geq n)$,
whereas
$\oo{P}_{\om,\beps}^o (\oo{R}_k < \infty, X_{\oo{S}_{k+1}} = z, 
\oo{S}_{k+1} = n) $ 
is measurable on
$\sigma (\om_x: x\cdot \ell \le z\cdot \ell -L |\ell|^2/|\ell|_1) \otimes 
\sigma(\eps_i, i < n)$.  Hence, by the
$\phi$-mixing property on cones of $P$, by the product structure of
$Q$ and by stationarity,
\begin{align}
\oo\PP^o (\oo{R}_{k+1}  < \infty)
& \le
\sum_{z\in\ZZ^d} \sum_{n\in\NN} \left[
 E_{P\otimes Q} 
\Bigl(\oo{P}_{\om,\beps}^o (\oo{R}_k < \infty, X_{\oo{S}_{k+1}} = z, 
\oo{S}_{k+1} = n) \Bigr)
\cdot E_{P\otimes Q} 
\Bigl(\oo{P}^o_{{\om, \beps}} (D' < \infty)
\Bigr) \right]  \nonumber \\
& + \phi( L )\sum_{z\in \ZZ^d} \sum_{n\in\NN} {E}_{P\otimes Q}
\Bigl( \oo{P}^o_{\om,\beps} 
(\oo{R}_k < \infty, X_{\oo{S}_{k+1}} = z, 
\oo{S}_{k+1} = n)
\Bigr)  \nonumber \\
& = \oo\PP^o (\oo{R}_{k}  < \infty)
(\oo\PP^o(D'<\infty)+\phi(L))  \nonumber \\ 
& \le 
(\oo\PP^o(D'<\infty)+\phi(L))^{k+1}
\label{taufinite2}
\end{align}
by induction. Choosing $L$ with $\phi(L)\leq \eta$, and 
using (\ref{boundD'}), we 
see that  $\oo{\PP}^o(K \geq k) \leq (1-\eta)^k$.
\qed
\end{proof}

 Consider now $\tau^{(L)}_1$ as a function of the path
$(X_n)_{n \ge 0}$ and set
\begin{equation} 
\label{tauk}
\tau^{(L)}_{k+1} = \tau^{(L)}_k (X_\bfcdot) + \tau^{(L)}_1 (X_{\tau^{(L)}_k +\bfcdot} - X_{\tau^{(L)}_k})\,,
\end{equation}
with $\tau^{(L)}_{k+1} = \infty$ on $\{\tau^{(L)}_k = \infty\}$. 

Under  (${\cal A} 3$),  $\tau^{(L)}_k$ is  $\oo{\PP}^o$-a.s. finite for all
$k$. Indeed, in view of the definition (\ref{tauk}),
\begin{align*}
\oo{\PP}^o(  \tau_1^{(L)}< \infty, \tau_2^{(L)} = \infty) &=
\sum_{z\in\ZZ^d} \sum_{n\in\NN} 
 E_{P\otimes Q}
\oo{P}_{\om,\beps}^o
( \tau_1^{(L)}=n, X_n=z, 
\tau_2^{(L)} = \infty)\\
& \leq 
\sum_{z\in\ZZ^d} \sum_{n\in\NN} 
 E_{P\otimes Q}
\oo{P}_{\om,\theta^n \beps}^z
( \tau_1^{(L)} = \infty)\\
 &= \sum_{z\in\ZZ^d} \sum_{n\in\NN} 
\oo{\PP}^z(  \tau_1^{(L)} = \infty) = 0\;,
\end{align*}
since all summands are equal to $\oo{\PP}^o(  \tau_1^{(L)} = \infty) = 0$.

Define
$$
\calH_1 = \sigma \left(\tau^{(L)}_1, X_0, \eps_0, X_1,\cdots ,
 \eps_{\tau^{(L)}_1-1},  X_{\tau^{(L)}_1},
\{\om (y, \cdot
); {\ell\cdot y <  \ell \cdot X_{\tau^{(L)}_1}}
-L|\ell|^2/|\ell|_1\}\right)  \,,
$$
$$
\calH_k =\sigma \Bigl(\tau^{(L)}_1 \ldots \tau^{(L)}_k, \quad X_0, \eps_0, X_1,
\cdots , \;  \eps_{\tau^{(L)}_k-1},X_{\tau^{(L)}_k},
\quad
\{ \om (y, \cdot ); {\ell \cdot y < \ell \cdot X_{\tau^{(L)}_k}}
-L|\ell|^2/|\ell|_1\} \Bigr)
\,.$$
Note that since
$\{D'=\infty\} = \{X_1, \ldots  X_{\tau^{(L)}_1} \in C(0, \ell, \zeta)\}
\bigcap  \{ \tau^{(L)}_1 < \infty \}$, see (\ref{stayincone}),
we have that
$$
\{D'=\infty\} \in \calH_1\;.
$$
Then, we have the following crucial lemma. Recall  the variational distance
 $\|\mu-\nu \|_{\rm var}=\sup\{ \mu(A)-\nu(A); A {\rm \ measurable}\}$ between two 
probability measures on the same space.
\begin{lemma}
\label{lem-"renewal"}
Assume (${\cal A} 1, 2, 3$), and  $\zeta \leq
\delta(\ell)/(3|\ell|)$. Set $\phi'(L)=2[\oo{\PP}^o (D'=\infty)-
\phi(L)]^{-1} \phi(L)$. (Here, and in 
the following, we consider $L$
large enough so that  $\phi(L)< \oo{\PP}^o (D'=\infty)$.)
Then, it holds a.s., 
$$
\|\oo{ \PP}^o
\Bigl( \{X_{\tau_k^{(L)}+n} - X_{\tau_k^{(L)}}\}_{n \ge 0} \in \cdot \; \vert
\calH_k \Bigr) - 
 \oo{ \PP}^o
\Bigl( \{X_n\}_{n \ge 0} \in \cdot \; \vert
D' = \infty \Bigr) \|_{\rm var} \leq \phi'(L)\;.
$$
\end{lemma}

\begin{proof}
We start with the case $k=1$. Let $A$ be a measurable subset of the path space,
and write for short ${\bf 1}_A={\bf 1}_{\{X_n-X_0\}_{n \geq 0} \in A}$.
Let $h \geq 0$ be a  $\calH_1$-measurable non-negative random variable.
Then for all $l, n \geq 1, x \in \ZZ^d$, there exists a random variable
$h_{x,l,n} \geq 0$, measurable with respect to 
$\sigma ( \{\om(y, \cdot);  y \cdot \ell  < x \cdot \ell -L |\ell|^2/|\ell|_1\},
\{X_i\})_{i \leq n}$ such that, on the event 
$\{ \tau_1^{(L)}=\oo{S}_l=n, X_{\oo{S}_l}=x\}$, it holds $h=h_{x,l,n}$. 
Recall that $\{K=l\}=\{
 \oo{S}_l < \infty, D' \circ \theta_{\oo{S}_l} = \infty\}$, and 
use the (weak) 
Markov property and shift invariance to write
\begin{align}
\nonumber
\oo{\EE}^o( h \; {\bf 1}_A \circ \theta_{\tau_1^{(L)}} ) &=
\sum_{l \geq 1} 
\oo{\EE}^o( h \; {\bf 1}_A 
\circ \theta_{\tau_1^{(L)}}  {\bf 1}_{K=l} )\\
\nonumber
&=
\sum_{l \geq 1, x \in \ZZ^d, n \geq 1}  E_{P\otimes Q}
\oo{E}_{\om,\beps}^o( 
h_{x, l, n}  \;{\bf 1}_A \circ \theta_n  
 {\bf 1}_{\oo{S}_l=n, X_n=x, D' \circ \theta_n}=\infty)\\
\nonumber
&=
\sum_{l \geq 1, x \in \ZZ^d, n \geq 1}  E_{P\otimes Q}\left[ 
\oo{E}_{\om,\beps}^o (
h_{x, l, n}  {\bf 1}_{\oo{S}_l=n, X_n=x}) \times
\oo{P}_{\om,\theta_n \beps}^x(A \bigcap \{D'=\infty\}) \right] \\
\nonumber
&=
\sum_{l \geq 1, x \in \ZZ^d, n \geq 1}
\oo{\EE}^o(h_{x, l, n}  {\bf 1}_{\oo{S}_l=n, X_n=x}) 
\oo{\PP}^x( A \bigcap \{D'=\infty\}) 
+ \rho_{A}\\
\label{432a}
&=
\oo{\PP}^o( A \bigcap \{D'=\infty\}) 
\sum_{l \geq 1, x \in \ZZ^d, n \geq 1}
\oo{\EE}^o(h_{x, l, n}  {\bf 1}_{\oo{S}_l=n, X_n=x}) + \rho_{A}\;.
\end{align}
The quantity $\rho_{A}$ is defined by the above equalities, i.e.,
$\rho_{A}=\sum_{l, x, n} {\rm Cov}_{P\otimes Q} 
(f_{x,l,n}, g_{x,n})$ with 
$f_{x,l,n}=
\oo{E}_{\om,\beps}^o (
h_{x, l, n}  {\bf 1}_{\oo{S}_l=n, X_n=x}), 
g_{x,n}=
\oo{P}_{\om,\theta_n \beps}^x(A \bigcap \{D'=\infty\})$. The point is that, 
from (\ref{fmix}) it holds for a non-negative $h$
\begin{align}
\label{432b}
|\rho_{A}| \leq \phi(L) 
\sum_{l \geq 1, x \in \ZZ^d, n \geq 1}  
\oo{\EE}^o( h_{x,l,n}  {\bf 1}_{\oo{S}_l=n, X_n=x})
\end{align}
uniformly in $A$. In particular for $A$ equal to 
the whole path space $({\ZZ^d})^{\NN}$, one gets
\begin{align}
\label{432c}
\sum_{l \geq 1, x \in \ZZ^d, n \geq 1}  
\oo{\EE}^o( h_{x,l,n}  {\bf 1}_{\oo{S}_l=n, X_n=x})
 \leq 
[\oo{\PP}^o (D'=\infty)-
\phi(L)]^{-1} \oo{\EE}^o( h)
\end{align}
as well as a formula for the sum in the left member above. 
Plugging this formula in 
 (\ref{432a}), one obtains 
\begin{align*}
|\oo{\EE}^o( h \:{\bf 1}_A \circ \theta_{\tau_1^{(L)}} ) - 
\oo{\EE}^o( h) \oo{\PP}^o (A | D'=\infty)| &= |\rho_{A} - \rho_{({\ZZ^d})^{\NN}}
\oo{\PP}^o (A | D'=\infty)|\\
& \leq 2 [\oo{\PP}^o (D'=\infty)-
\phi(L)]^{-1} \phi(L) \oo{\EE}^o( h)
\end{align*}
where the second inequality follows from
(\ref{432b}), (\ref{432c}).
Since $h$ is arbitrary, we have
$$
\vert  \oo{ \PP}^o
\Bigl( \{X_{\tau_1^{(L)}+n} - X_{\tau_1^{(L)}}\}_{n \ge 0} \in A \vert
\calH_1 \Bigr) - 
 \oo{ \PP}^o
\Bigl( \{X_n\}_{n \ge 0} \in A \vert
D' = \infty \Bigr) \vert \leq \phi'(L)\,.
$$
a.s., for all $A$'s. But there are only countably many cylinders in the path 
space, 
so we can find a subset of $\Omega$ of $P$-measure one, where the previous 
inequality holds simultaneously
for  all measurable $A$. We have shown the lemma for $k=1$.

The case of a general $k \geq 1$ follows similarly 
from the above computations, and the  
definition (\ref{tauk}).
\qed
\end{proof}

\section{Law of large numbers}
\label{sec:LLN}
\setcounter{equation}{0}
Throughout this section we assume $({\cal A}1,2)$ for some $\ell$
with integer coordinates, and we assume also that the conclusions of
Lemma \ref{lem-"renewal"} hold.
For $L \in |\ell|_1 \NN^*$ we  define  $\tau^{(L)}_0=0$, and for $k\geq 1$,
\begin{equation} 
\label{bartau}
{\ot}_k^{(L)}= \kappa^L \left(\tau^{(L)}_k - \tau^{(L)}_{k-1}\right)\,,
\quad\quad
{\oX}_k^{(L)}= \kappa^L \left(X_{\tau^{(L)}_k} - X_{\tau^{(L)}_{k-1}}\right)\,.
\end{equation}

The following uniform integrability condition is instrumental in our derivation:
\begin{assumption}
\label{ass-C}
\mbox{}
\begin{itemize}
\item[(${\cal A} 5$)]
There exist  an $\alpha>1$ and
$M=M(L)$ such that 
$\phi'(L)^{1/\alpha'} M(L)^{1/\alpha}
\underset{L\to\infty }{\longrightarrow}
0$ 
(with $1/\alpha'=1-1/\alpha$), and
\begin{equation}
\label{stronger}
P\Big(\;\oo \EE^o\big((\ot_1^{(L)})^\alpha\;|\;D'=\infty, {\cal F}_0^L\big)>M\;\Big)=0\,,
\end{equation}
where ${\cal F}_0^L=\sigma(\omega(y,\cdot): \ell\cdot y<-L)$.
\end{itemize}
\end{assumption}
We define
\begin{equation}
\label{integ}
\beta_L:=\oo \EE^o(\ot_1^{(L)}\;|\;D'=\infty)<\infty\,,
\end{equation}
and 
\begin{equation}
\label{integ2}
\gamma_L:=\oo \EE^o(\oX_1^{(L)}\;|\;D'=\infty) \in \reals^d \,,
\end{equation}
where
the moments $\beta_L$ and $\gamma_L$ are finite
due to (\ref{stronger}). Further, we note that $\beta_L\geq 1$ 
for all $L$.

{}From Lemma \ref{lem-"renewal"}, we have a.s., that
for $k\geq 2$,
\begin{equation}
\label{ff}
\| \oo{ \PP}^o\left(({\ot}_k^{(L)}, {\oX}_k^{(L)})\in \cdot\; |\;
{\cal H}_{k-1}\right)- \mu^{(L)}( \cdot )
\|_{\rm var} \leq  \phi'(L)\,,
\end{equation}
where $\mu^{(L)}$ is defined by
$$\mu^{(L)}(A\times B)=
\oo{\PP}^o\left({\ot}_1^{(L)}\in A, {\oX}_1^{(L)}\in B\,|\, D'=\infty
\right)\,,$$
for any sets $A \subset \kappa^L \NN^*, B \subset \kappa^L \ZZ^d$.
This  will allow us to implement a coupling procedure. 
We recall the following splitting representation: if ${\overline X},{\tilde X}$ 
are random variables
of laws ${\overline P}, {\tilde P}$ such that 
$\|  {\overline P} - {\tilde P}\|_{\rm var}\leq
a$ then one may find, on an enlarged probability space, 
independent random variables $Y,\Delta,
Z,{\tilde Z}$ where  $\Delta $ is Bernoulli distributed on $\{0,1\}$ with
parameter $a$, and 
$$
{\overline X} =(1-\Delta) Y + \Delta Z\;,\;\;
{\tilde X} =(1-\Delta) Y + \Delta {\tilde Z}\,,$$
(see e.g. \cite[Appendix A.1]{barbour} for the proof);
in particular, 
$${\overline X}=(1- \Delta) {\tilde X}+ \Delta Z\;,
\quad
|\Delta Z| \leq | {\overline X}|, \;
|\Delta {\tilde Z}| \leq | {\tilde X}| \,.
$$
Thus, due to  (\ref{ff}) --see for similar constructions 
\cite{berbee} or \cite[Chapter 3]{thorisson}--,
we can enlarge our probability space where is defined the 
sequence $\{(\ot_i^{(L)},\oX_i^{(L)})\}_{i\geq 1}$ in order to support
also:

\begin{itemize}
\item a  sequence 
 $\{(\tilde \tau_i^{(L}, \tilde X_i^{(L)}, \Delta_i^{(L)})\}_{i\geq 1}$
of i.i.d. random vectors (with values 
in $\kappa^L \NN^* \times \kappa^L \ZZ^d \times \{0,1\}$)
such that  
$\{(\tilde \tau_1^{(L)}, \tilde X_1^{(L)})\}$ is distributed
according to $\mu^{(L)}$ while
$\Delta_1^{(L)}\in \{0,1\}$
is such that $P(\Delta_1^{(L)}=1)=\phi'(L)$,

\item and another sequence 
 $\{(Z_i^{(L)},Y_i^{(L)})\}_{i\geq 1}$ such that
$$(\ot_i^{(L)},\oX_i^{(L)})=
(1-\Delta_i^{(L)})(\tilde\tau_i^{(L)},\tilde X_i^{(L)})+
\Delta_i^{(L)} (Z_i^{(L)},Y_i^{(L)})\,,$$
and such that, with
$${\cal G}_i=\sigma(\{\tilde \tau_j^{(L}\}_{j\leq i-1},
\{\tilde X_j^{(L)})\}_{j\leq i-1},
\{\Delta_j^{(L)}\}_{j\leq i-1})\,,$$
it holds that $\Delta_i^{(L)}$ is independent of ${\cal G}_i$ and of
$(Z_i^{(L)}, Y_i^{(L)})$.
\end{itemize}

 The joint law of the variables  
$\{(Z_i^{(L)},Y_i^{(L)})\}_{i\geq 1}$
is complicated,  but
 it holds that $|Y_i^{(L)}|\leq Z_i^{(L)}$
while, due to (\ref{stronger}) and since $|\Delta_i^{(L)} Z_i^{(L)} |
\leq \tau_i^{(L}$, 
it holds almost surely that
\begin{equation}
\label{test1}
\oo{\EE}^o[ ( \Delta_i^{(L)} Z_i^{(L)})^\alpha\,|\, {\cal G}_i]=
\phi'(L)\oo{\EE}^o[ ( Z_i^{(L)})^\alpha\,|\, {\cal G}_i]\leq M(L)\,.
\end{equation}   

We next have the 
\begin{lemma}
\label{mainconv}
Assume the integrability condition (\ref{stronger}). Then,
there exists a sequence $\eta_L
\underset{L\to\infty }{\longrightarrow} 0$ such that
\begin{equation}
\label{timeconv}
\limsup_{n\to\infty}
\left|\frac{1}{n} \sum_{i=1}^n \ot_i^{(L)}-\beta_L\right|<\eta_L\,,\quad
\oo\PP^o-a.s.,
\end{equation}
and
\begin{equation}
\label{spaceconv}
\limsup_{n\to\infty}
\left|\frac{1}{n} \sum_{i=1}^n \oX_i^{(L)}-\gamma_L\right|<\eta_L\,,\quad
\oo\PP^o-a.s.,
\end{equation}

\end{lemma}
\noindent
{\bf Proof of Lemma \ref{mainconv}}
We prove (\ref{timeconv}), the proof of (\ref{spaceconv}) being
similar. Simply write 
$$
\frac{1}{n} \sum_{i=1}^n \ot_i^{(L)} = \frac{1}{n} \sum_{i=1}^n
\tilde\tau_i^{(L)}
-\frac{1}{n} \sum_{i=1}^n \Delta_i^{(L)} \tilde\tau_i^{(L)}
+\frac{1}{n} \sum_{i=1}^n \Delta_i^{(L)}  Z_i^{(L)}\;.
$$
Note first that by independence,
$$ \frac1n \sum_{i=1}^n  \tilde\tau_i^{(L)}
\underset{n\to\infty }{\longrightarrow}
\beta_L\,,\quad \oo \PP^o-a.s.,$$
while 
\begin{equation}
\label{yom}
 |\frac1n \sum_{i=1}^n\Delta_i^{(L)} \tilde\tau_i^{(L)} |
\leq 
\left(
\frac1n \sum_{i=1}^n(\Delta_i^{(L)})^{\alpha'}\right)^{1/\alpha'}
\left(\frac1n \sum_{i=1}^n(\tilde \tau_i^{(L)})^\alpha\right)^{1/\alpha}
\end{equation}
and hence 
$$\limsup_{n\to\infty} 
|\frac1n \sum_{i=1}^n\Delta_i^{(L)} \tilde\tau_i^{(L)} |
\leq \phi'(L)^{1/\alpha'} M(L)^{1/\alpha}, \quad \oo \PP^o-a.s.$$
We next consider the term involving $Z_i^{(L)}$. 
Set $\oo Z_i^{(L)}:=\oo\EE^o(Z_i^{(L)}|{\cal G}_i)$, and note that
$M_n:=\sum_{i=1}^n \Delta_i^{(L)}
 (Z_i^{(L)}-\oo Z_i^{(L)})/i$ is a 
zero mean 
martingale with respect to the filtration ${\cal G}_i$. 
By the Burkholder-Gundy
maximal inequality  \cite[14.18]{williams}, for $\gamma=\alpha \wedge 2$,
$$E |\sup_n  M_n|^\gamma\leq C_\gamma E\left(\sum_i
 \frac{(\Delta_i^{(L)}
 (Z_i^{(L)}-\oo Z_i^{(L)}))^2}{i^2}\right)^{\gamma/2}  
\!\!\!\!\!\!\!\!\!
\leq C_\gamma\sum_i E\left(\frac{(\Delta_i^{(L)}  
 (Z_i^{(L)}-\oo Z_i^{(L)}))^\gamma}{i^\gamma}\right)
\leq C_\gamma'\,,$$
for some constants $C_\gamma,C_\gamma'$. Hence, $M_n$ converges
$P$-a.s. to an integrable random variable, 
and by the Kronecker lemma \cite[12.7]{williams}, it holds that $n^{-1}\sum_{i}
\Delta_i^{(L)}
 (Z_i^{(L)}-\oo Z_i^{(L)})\to 0$, almost surely.
On the other hand, 
there is nothing to prove if $\phi'(L)=0$ while, if
$\phi'(L)>0$ then
$$|\oo Z_i^{(L)} |\leq \left(\oo \EE^o(|Z_i^{(L)}|^\alpha\,|
\,{\cal G}_i)
\right)^{1/\alpha}
\leq \left(\frac{M(L)}{\phi'(L)}\right)^{1/\alpha}\,$$
by (\ref{stronger}), and hence
$$ |\frac1n \sum_{i=1}^n \oo Z_i^{(L)} \Delta_i^{(L)}|
\leq \left(\frac{M(L)}{\phi'(L)}\right)^{1/\alpha}
\frac1n \sum_{i=1}^n  \Delta_i^{(L)}
\underset{n\to\infty }{\longrightarrow}
M(L)^{1/\alpha} \phi'(L)^{1/\alpha'}\,,\;
 \PP^o-a.s.\,,$$
yielding (\ref{timeconv}) by choosing 
$\eta_L=2  M(L)^{1/\alpha} \phi'(L)^{1/\alpha'}
\,.$
\qed

Using that 
$\tilde \tau_i^{(L)} \geq \kappa^L$ and that $\beta_L\geq 1$,
we conclude from Lemma \ref{mainconv} that for all $L$ large enough,
$$\limsup_{n\to\infty}\left|
\frac{ \frac1n \sum_{i=1}^n {\oo X}_i^{(L)}}
{\frac1n \sum_{i=1}^n {\oo \tau}_i^{(L)}}-\frac{\gamma_L}{\beta_L}\right|
\leq 3 \eta_L\,, \quad \oo \PP^o-a.s.,$$
from which one deduces by standard arguments that
$$\limsup_{n\to\infty} \left|\frac{X_n}{n} -\frac{\gamma_L}{\beta_L}\right|
\leq 4\eta_L\,,\quad \oo\PP^o-a.s.\,.$$
Thus, we conclude both the existence of the limit $v:=
\lim_{L\to\infty} \gamma_L/\beta_L$ and the $\oo\PP^o$ convergence 
of $X_n/n$ to it.
Summarizing, we have proved the following
\begin{theorem}
\label{theo-LLN}
Assume the conclusion of Lemma \ref{lem-"renewal"} and the integrability
condition (\ref{stronger}). Then, there exists a deterministic 
vector $v$ with $v \cdot \ell >0$ such that
$$\lim_{n\to\infty}
\frac{  X_n}{n}=v\,,\quad \PP^o-a.s..$$
\end{theorem}

{\bf Proof
of Theorem \ref{theo-LLN}}
The only statement left  to be shown is $v \cdot \ell >0$. Actually
we show that for all $L \in |\ell|_1 \NN^*$, 
we have $\gamma_L = v \beta_L$, from 
which the desired claim follows, since $\gamma_L \cdot \ell \geq \kappa^L$ and 
$\beta_L < \infty$. Let us fix $L$. 
We already know that $ X_n/n \to v$, $\PP^o$- a.s.,
and since 
$\tau_n^{(L)} \to \infty$, we have
$$
 \frac{  X_{\tau_n^{(L)}}}{\tau_n^{(L)}} -v \; = 
\left( \frac{  X_{\tau_n^{(L)}}}{n} - v \frac{\tau_n^{(L)}}{n} \right)
\frac{n}{  \tau_n^{(L)}} 
\underset{n\to\infty }{\longrightarrow}
 0\,, \quad \PP^o-a.s. 
$$
By (\ref{stronger}) $\oo \EE^o \tau_n^{(L)}/n \leq M(L)
 \kappa^{-L}$, so that
$\tau_n^{(L)}/n$ is bounded in probability, yielding
 \begin{equation}
\label{pesah1} 
\lim_{n \to \infty} \left( \frac{  X_{\tau_n^{(L)}}}{n} - v 
\frac{\tau_n^{(L)}}{n}
\right) = 0\,, \quad \PP^o-a.s. 
\end{equation}
Again from (\ref{stronger}) and Jensen's inequality, we have 
 \begin{equation*}
\oo \EE^o \left( \frac1n { [\tau_{n+1}^{(L)}-\tau_1^{(L)}]} \right)^{\alpha } 
 \leq 
\oo \EE^o \frac1n {\sum_{i=1}^n  \left( \kappa^{-L}
\bar \tau_{i+1}^{(L)} \right)^{\alpha }} \leq M(L)\kappa^{-\alpha L}\;.
\end{equation*} 
So the sequence of
random variables $\{[\tau_n^{(L)}-\tau_1^{(L)}]/n\}_{n \geq 1}$ is 
uniformly integrable,
and the similar conclusion holds for $\{[X_{\tau_n^{(L)}} -X_{\tau_1^{(L)}}]
/n\}_{n \geq 1}$.
Therefore, taking the expectation in (\ref{pesah1}),
we obtain
$$
\kappa^{-L}( \gamma_L - v \beta_L) =0
$$
for arbitrary $L$.
\qed

\noindent
{\bf Proof    
of Corollary \ref{cor-nonnest}}. 
Throughout, we write $\oo{P}_\om^z$ for $E_{Q}
\oo{P}_{\om,\beps}^z$, and assume w.l.o.g. that
$\kappa$ is taken small enough such that $\delta(\ell)>2\kappa$. 
(Note that $\oo{P}_\om^z$, in contrast to  ${P}_\om^z$,
 allows considering the $\beps$ random sequence.)
Recall that from Lemma
\ref{D'>0},  
$$
q=\inf_{\om   } \oo{P}_{\om}^o (D'  = \infty )=
\inf_\om P_\om^o(D'=\infty)\in (0,1]\,.
$$
For $a>0$, set
$$\phi_0^L(a)=\sup_\om \oo{E}_\om^o(\exp a \kappa^L\oo{S}_1)\,,\quad
\phi_1^L(a)=\sup_\om \oo{E}_\om^o(\exp a \kappa^L\oo{S}_1 \,|\, 
D'<\infty)\,.$$
We derive below estimates on $\phi_0^L(a),\phi_1^L(a)$ 
and in particular
show that under a non-nestling assumption, these functions are 
bounded
(uniformly in $L$) 
for $a$ small enough, and hence can be made arbitrarily close
(again, uniformly in $L$)
to $1$ by reducing $a>0$. Assuming that, 
we get from the Markov property and independence 
of the $\beps$ sequence, 
\begin{eqnarray}
\label{pesah2}
\oo{E}_\om^o(\exp \{a \ot_1^{(L)})\}&=&
\sum_{k=1}^\infty 
\oo{E}_\om^o(\exp \{a \ot_1^{(L)}\} \, {\bf 1}_{K=k})\nonumber\\
&=& \sum_{z\in \ZZ^d}
 \oo{E}_\om^o( \exp \{a\kappa^L\oo{S}_1 \} {\bf 1}_{X_{\oo{S}_1}=z})
\cdot \oo{P}_\om^z({D'=\infty})\nonumber\\
&& \!\!\!\!\!\!\!\!\!\!\!
\!\!\!\!\!\!\!\!\!\!\!
 \!\!\!\!\!\!\!\!\!\!\!
\!\!\!\!\!\!\!\!\!\!\!
+ \sum_{z,z_1\in \ZZ^d}
 \oo{E}_\om^o( \exp \{ a\kappa^L\oo{S}_1\}  {\bf 1}_{X_{\oo{S}_1}=z})
\cdot \oo{P}_\om^z({D'\!<\!\infty})
 \oo{E}_\om^z( \exp \{a\kappa^L\oo{S}_1 \}
 {\bf 1}_{X_{\oo{S}_1}=z_1} | D'\!<\!
\infty)
\cdot \oo{P}_\om^{z_1}({D'\!=\!\infty})+ \ldots\nonumber\\
&\leq &
\sup_\om \oo{P}_\om^o({D'=\infty})
\left[\phi_0^L(a)+\phi_0^L(a) \phi_1^L(a)(1-\inf_\om 
\oo{P}_\om^o({D'=\infty}))+\ldots\right]
\nonumber\\
&\leq& \phi_0^L(a)\left[\sum_{k=0}^\infty \phi_1^L(a)^k
(1-q)^k
\right]\nonumber\\
&= &\frac{\phi_0^L(a)}{(1-(1-q)\phi_1^L(a))}:=g(a,L)\, 
\end{eqnarray}
where $\sup_{L} g(a,L)<\infty$ for small enough, positive $a$.
Thus,
$$
\EE^o((\ot_1^{(L)})^\alpha|D'=\infty, {\cal F}_0^L) \leq
\frac{ \sup_{\om} \oo{E}_{\om}^o
[(\ot_1^{(L)})^\alpha]}
{ \inf_{\om } \oo{P}_{\om}^o
(D'=\infty) } \leq
\frac{ \C}{ q a^\alpha} \sup_{\om, L }
 \oo{E}_{\om}^o  e^{a \ot_1^{(L)}}\;,
$$
which, together 
with (\ref{pesah2}) and choosing $a>0$ small enough,
 yields (\ref{stronger}) with any $\alpha >1$ and some $M$
not depending on $L$. Hence,
Theorem \ref{theo-LLN} applies for 
any rate $\phi \to 0$.

We thus turn to the proof of the claimed (uniform in $L$)
 finiteness of $\phi_0^L(a)$ 
and $\phi_1^L(a)$ for $a$ small enough. Since 
$\phi_1^L(a)\leq \frac{\phi_0^L(a)}{1-q}\,$
it clearly suffices to consider $\phi_0^L(a)$.
Let us denote here by 
$T_m  (m=1,2\ldots)$ the hitting time of the half-space
$\{ x \cdot \ell \geq m L |\ell|^2/ |\ell|_1 \}$, limited by the
hyperplanes
through points $m ( L/ |\ell|_1)\ell$ and orthogonal to $\ell$.
Time $T_m$ is called \textit{L-successful\/} if
$$(\eps_{T_m+1}, \eps_{T_m+2},\cdots \eps_{T_m+L})= 
\bar \eps^{(L)}\;.$$
We denote by  $I=\inf \{ m \geq 1; T_m\; {\rm is \ L-successful}\}$, and we 
note 
that, by definition,
$ \oo{S}_{1} \leq T_I + L$, and that $I$ is geometrically 
distributed on $\NN^*$ with failure probability 
$\kappa^L$. Let 
$$\psi_0^L(a)=\sup_\om \oo{E}_\om^o(\exp a \kappa^L (T_1+L))\,,\quad
\psi_1^L(a)=\sup_\om \oo{E}_\om^o(\{\exp a \kappa^L T_1 \} \, {\bf 1}_{\{
(\beps_1, \cdots  \beps_L) \neq \bar \eps^{(L)}\}})
\,.$$
Similar to (\ref{pesah2}), 
\begin{eqnarray*}
\oo{E}_\om^o(\exp \{a \kappa^L \oo{S}_1\}&\leq &
\sum_{m=1}^\infty 
\oo{E}_\om^o(\exp \{a \kappa^L (T_m+L)
\} \, {\bf 1}_{I=m})\\
&=& 
\oo{E}_\om^o(\exp \{a \kappa^L (T_1+L)\}) \kappa^L
+ \\
\\
&& 
+ \sum_{z\in \ZZ^d} 
\oo{E}_\om^o(\exp \{a \kappa^L (T_1+L)\} {\bf 1}_{X_{T_1}=z})
 \oo{E}_\om^o(\exp \{a \kappa^L (T_2-T_1)\} {\bf 1}_{I>1})\kappa^L+\ldots\\
& \leq & 
\kappa^L \psi_0^L(a)\left[\sum_{k=0}^\infty \psi_1^L(a)^k
\right]= 
\frac{\kappa^L \psi_0^L(a)}{(1-\psi_1^L(a))_+}\,,
\end{eqnarray*}
with $(r)_+= \max\{r, 0\}$. 
But, if $a \leq \lambda_0 \delta(\ell) \kappa^{-|\ell|_1}$,
\begin{eqnarray*}
\psi_1^L(a) & \leq &
(1- \kappa^L)
\sup_{\om, \beps} \oo{E}_{\om, \beps}^o(\{\exp a \kappa^L T_1 \}\\
 & \leq &(1- \kappa^L)
\exp\{ 3 a \kappa^L/ \delta(\ell)\}\;,
\end{eqnarray*}
from (\ref{expexit-q}). Hence, we can choose $a>0$ small enough so that 
$\psi_1^L(a)< \infty$ for all $L$ and such that
$$\sup_L \frac{\kappa^L }{(1-\psi_1^L(a))_+}< \infty\,, \quad
\sup_L \psi_0^L(a) < \infty.$$
 This implies that
$\sup_L \phi_0^L(a) < \infty$ for $a>0$ small.
\qed

\section{Mixing}
\label{sec:mix}
\setcounter{equation}{0}

Here are the main examples of distributions $P$ of the environment 
field which are  $\phi$-mixing on cones.

\begin{definition}

\begin{enumerate}
 \item A random field $P$ is $\phi$-mixing if 
there exists a function
$\phi (r) \underset{r\to\infty}{\to} 0$ such that any two
$r$-separated events $A, B$ with $P(A) > 0$,
$$
\left| \frac{P(A \cap B)}{P(A)} - P(B) \right| \le \phi (r)\,.
$$
 \item Let $k \geq 1$, and let $\partial \La^k=\{ z \in \lambda^c; 
\dist(z,
   \La) \leq k\}$
be the $k$-boundary of $\La \subset \ZZ^d$. ($\dist$ and $|.|$ both
denote the Euclidean distance).
A random field $P$ is {\bf $k$-Markov}  if 
there exists a family $\pi$ of transition kernels --- called {\it
specification} --- $\pi_{\La}=\pi_{\La}( \prod_{y \in \La} d\om_y |\;
{\calF}_{\partial \La})$ for finite  $ \La \subset \ZZ^d $
such that
\begin{equation} 
\label{markov}
P \big( (\om_x)_{x \in \La} = \cdot \; |\; {\calF}_{ \La^c}
\big) =
 \pi_{\La} \big(\;\cdot \;| \; {\calF}_{\partial \La}
\big), 
\quad P-\mbox{-a.s.}
\end{equation}
In addition, a $k$-Markov field 
$P$  is called {\bf weak-mixing} if there
exist constants
$\gamma >0$, $ C< \infty$ such that for all finite subsets
$V \subset \La \subset \ZZ^d$,
\beq{weakmix}
  \sup \big\{ \| 
\pi_{\La}({\cdot\mid \om }) - \pi_{\La}({\cdot \mid \om'})
\|_{V};
  \quad \om, \om' \in {\Si}^{\La^c} \big\} \leq C \sum_{y \in V\!,\, z \in
\partial
  \!  \La^k} \exp (-\gamma |z-y|) \;,
\end{equation}
with $\|.\|_V=\|.\|_{{\rm var},V}$ the variational norm on $V$, $\|\mu-\nu\|_V=\sup\{
\mu(A)-\nu(A); A \subset \si((\om_x)_{x \in V})\}$.
\end{enumerate}
\end{definition}

These notions of mixing are different and both of practical interest.
Refer to \cite{doukhan} for the first one. The second one describes
environments produced by a Gibbsian particle system at equilibrium
in the uniqueness regime \cite{dobrushinshlosman, martinelli}.

\begin{proposition}
\label{d51}
Assume $P$ is stationary and ergodic.
If $P$ is $\phi$-mixing, 
 then $P$ is $\phi$-mixing on cones, i.e.
Assumption   (${\cal A} 1$) is satisfied. 
When $P$ is weak-mixing $k$-Markov of constant $\gamma$, then
Assumption   (${\cal A} 1$) is satisfied  with
the function $\phi(r)=\C(\zeta)  e^{-\gamma' r}$, and 
$\gamma'=\gamma/ \sqrt{2}$.
\end{proposition}

\begin{proof}
\begin{enumerate}
 \item In the $\phi$-mixing case, the statement  directly  follows
   from the definition. We can even take $\zeta=0$, i.e. we can
   replace cones by hyperplanes as in \cite{sznitmanzerner}, all
   through the paper.
 \item We assume now that $P$ is a weak-mixing $k$-Markov field. 
Fix some 
$\zeta>0$ and $\ell \in \reals^d \setminus \{0\}$.
For $m, M, N >0$, define the truncated cone $V$
$$
V_m=V=C(r \ell, \ell, \zeta) \bigcap \{ y \in \reals^d ; 
(y-r\ell) \cdot \ell \leq m |\ell|^2 \}\;,
$$ 
and  the cylinder $\Lambda$ 
$$
\La=\{ y \in \reals^d ; 
0 \leq y \cdot \ell \leq (r\!+\!m\!+\!M)|\ell|^2  , \vert y-(y\cdot \ell) 
\ell/|\ell|^2 \vert
\leq N|\ell|\! +\! m|\ell| \tan( \cos^{-1}(\zeta))\}\;,
$$
depicted in Figure \ref{fig-1a}.

\begin{figure}[h]
\begin{center}
\begin{picture}(0,0)%
\includegraphics{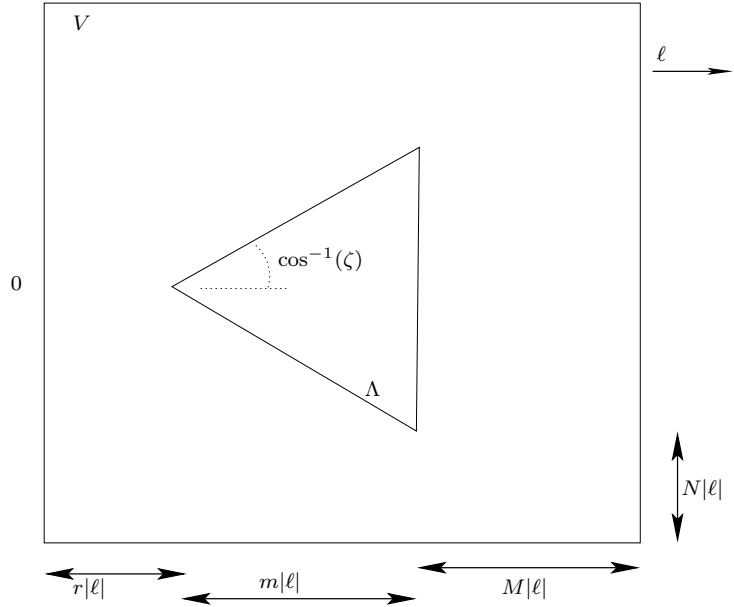}%
\end{picture}%
\setlength{\unitlength}{2052sp}%
\begingroup\makeatletter\ifx\SetFigFont\undefined%
\gdef\SetFigFont#1#2#3#4#5{%
  \reset@font\fontsize{#1}{#2pt}%
  \fontfamily{#3}\fontseries{#4}\fontshape{#5}%
  \selectfont}%
\fi\endgroup%
\begin{picture}(8722,7294)(1951,-7793)
\put(2701,-811){\makebox(0,0)[lb]{\smash{\SetFigFont{8}{9.6}{\familydefault}{\mddefault}{\updefault}$V$}}}
\put(6226,-5236){\makebox(0,0)[lb]{\smash{\SetFigFont{8}{9.6}{\familydefault}{\mddefault}{\updefault}$\Lambda$}}}
\put(1951,-3961){\makebox(0,0)[lb]{\smash{\SetFigFont{8}{9.6}{\familydefault}{\mddefault}{\updefault}0}}}
\put(2701,-7636){\makebox(0,0)[lb]{\smash{\SetFigFont{8}{9.6}{\familydefault}{\mddefault}{\updefault}$r |\ell|$}}}
\put(4951,-7561){\makebox(0,0)[lb]{\smash{\SetFigFont{8}{9.6}{\familydefault}{\mddefault}{\updefault}$m |\ell|$}}}
\put(7876,-7636){\makebox(0,0)[lb]{\smash{\SetFigFont{8}{9.6}{\familydefault}{\mddefault}{\updefault}$M |\ell|$}}}
\put(10051,-6436){\makebox(0,0)[lb]{\smash{\SetFigFont{8}{9.6}{\familydefault}{\mddefault}{\updefault}$N |\ell|$}}}
\put(9751,-1186){\makebox(0,0)[lb]{\smash{\SetFigFont{8}{9.6}{\familydefault}{\mddefault}{\updefault}$\ell$}}}
\put(5176,-3661){\makebox(0,0)[lb]{\smash{\SetFigFont{8}{9.6}{\familydefault}{\mddefault}{\updefault}$\cos^{-1}(\zeta)$}}}
\end{picture}
\caption{Cone definition}
\label{fig-1a}
\end{center}
\end{figure}

Then, $V \subset \Lambda$, and we split the sum over $z$ in
(\ref{weakmix}) into three terms, according to $
\partial \Lambda^k = 
 (\partial_- \Lambda) \bigcup  (\partial_+ \Lambda) \bigcup
 (\partial_O \Lambda) $ with
\beaa
 \partial_- \Lambda=  \partial \Lambda^k \bigcap \{ z \cdot \ell \leq
 0\}\;,\;
 \partial_+ \Lambda=  \partial \Lambda^k \bigcap \{ z \cdot \ell \geq
 (r+m+M)|\ell|^2\}\\ 
 \partial_O \Lambda =  \partial \Lambda^k \bigcap \{ \vert z-(z\cdot \ell) 
\ell/|\ell|^2 \vert \geq N|\ell|\! +\! m|\ell| \tan( \cos^{-1}(\zeta)
\}\;.
\eeaa
Since  $\zeta >0$, the number of points in the cone $C(r \ell, \ell,
\zeta)$ at distance $l$ from the hyperplane $z \cdot \ell =0$ grows
linearly in $l$. Using $\sqrt{ {a}^2+b^2} \geq (a+b)/\sqrt 2$ and 
the previous remark, we find that 
\beaa
\sum_{y \in V, z \in \partial_- \Lambda} \exp (- \gamma |z-y|) &\leq&
\sum_{y \in V} \C \exp (- \gamma' y\cdot \ell )\\ &\leq& 
\C  \exp (- \gamma'r |\ell|)
\eeaa
with 
$\C$ some constant (depending on
$k, \zeta$) which may
change from line to line. Similarly,
$$
\sum_{ y \in V, z \in \partial_O \Lambda}  \exp (- \gamma |z-y|) \leq
\C (N+m)^{d-2} \exp (- \gamma' N |\ell|) \;,
$$
and 
$$
\sum_{ y \in V, z \in \partial_+ \Lambda}  \exp (- \gamma |z-y|) \leq
\C m^{d-1} \exp (- \gamma' M |\ell|) \;.
$$
Letting now $M \to \infty$ and then $N \to \infty$ in (\ref{weakmix}), 
we get for 
all $A \in \si((\om)_x, x \cdot \ell \leq 0)$, and all 
$B \in \si((\om)_x, x\in V_m)$,  from (\ref{weakmix})
$$
|P(B|A)-P(B)| \leq \C  \exp (- \gamma' r |\ell|)\;,
$$
with arbitrary $m$.
But for any $B \in \si((\om)_x, x\in C(r \ell,
\ell, \zeta))$, there exists a sequence $B_m \in \si((\om)_x, x\in V_m)$
with $\lim_m P(B \Delta B_m)=0$. 
Therefore such a $B$ satisfies also the previous estimate, and
 $P$ satisfies 
(${\cal A} 1$) with exponentially vanishing $\phi$.
\qed
\end{enumerate}
\end{proof}

\section{A nestling example}
\label{sec-A5}
\setcounter{equation}{0}
For simplicity, we  work here with $\ell=e_1$. The example below can
be modified to work for any $\ell$, at the expense of more cumbersome
notations. Note that in this case, $({\cal A} 2)$ can be rephrased as
the directional ellipticity condition:

\noindent
$({\cal A} 2')$: There exists a $\kappa>0$ such that
$P(\om(0,e_1)>\kappa)=1$.

Our goal is to provide a family of
 \textit{examples\/} of environments which are nestling and 
to which the results of this paper apply. The examples can be 
considered as a perturbation of the environment with
$\om(x,x+e_1)=1$. Alternatively, they may also be 
considered as a perturbation of i.i.d. environments satisfying 
Kalikow's condition, perturbed by a slight dependence. 
(Indeed, it will appear from the proof
that $\delta_0 \searrow 0$ as $\gamma' \to + \infty$ in the statement
below.)

We claim the following:
\begin{theorem}
\label{theo-A5}
Assume $P$ satisfies $({\cal A}1)$ with $\phi(L)\leq \C(\zeta)
 e^{-\gamma' L}$ and
$({\cal A}2')$. Then there exists
a $\delta_0=\delta_0(\kappa,\gamma',d)<1$ such that 
\begin{center}
if
$P$ satisfies $({\cal A}3)$ with $\delta(e_1)>\delta_0$, then
$P$ satisfies $({\cal A}5)$.
\end{center}
\end{theorem}
A class of explicit examples satisfying the conditions of
Theorem \ref{theo-A5} is provided at the end of this section, 
following the:

\noindent
{\bf Proof of Theorem \ref{theo-A5}:} It is useful to consider the
marked point process of fresh times and  fresh points. Formally,
fresh times are those times when the random walk achieves a new record value
in the $\ell$ direction.
\begin{definition}
\label{def-1}
A point $x\in \ZZ^d$ is called a
 \textit{fresh point\/} for the RWRE $(X_n)$,
and a time $s$ is called a \textit{fresh time\/} for $(X_n)$, if
$$
\{ X_n \cdot \ell < x \cdot \ell, \quad n < s\}
\cap \{ X_s = x\}
\,.
$$
\end{definition}
We label these random couples $(s,x)$ according to increasing times,
$0=s_0<s_1< \ldots, x_n=X_{s_n} (n \geq 0)$. 
For transient walks in the direction $\ell$, 
there are infinitely many fresh times 
$\{s_i\}_{i \geq 0}$ and fresh points $\{x_i\}_{i \geq 0}$, 
and, in the present case $\ell = e_1$, 
it holds $x_{i+1}\cdot \ell = x_i \cdot \ell + 1$.

Like in the proof of Corollary \ref{cor-nonnest},
call a fresh time $s$ \textit{L-successful\/} if
$\eps_{s+1} = \eps_{s+2} = \cdots = \eps_{s+L} = e_1$.
Note that all $\oo{S}_k$ are L-successful fresh times, and that
an L-successful fresh time $s$ leads to an 
L-regeneration time (more accurately, an ``approximate L-regeneration time")
$s+L$ if $\theta_{s+L} D'=\infty$.

Define $F \geq 0$ by  $L+s_F=\tau_1^{(L)}$, so that it holds 
that
$$
X_{\tau_1^{ (L)}} \cdot \ell= F+L\;.
$$  
As a general feature, fresh points have much nicer tail properties
than  fresh times.
The following summarizes some properties
of fresh points and regeneration positions, and does not require 
any additional assumptions. It is slightly stronger than what we need
in the sequel.
\begin{lemma}
\label{lem-1}
Assume $({\cal A}1,2,3)$. Then there exist deterministic constants
$\zeta_0>0$ and 
$\la_2 = \la_2 (\delta(\ell), \kappa, d)$ and a function
$Q(\la)$ (depending on $\delta(\ell), \kappa$)
 such that for any $\zeta<\zeta_0$,
$\la<\la_2$ and all $L> L_0 (\la)$,
$$
\oo{\EE}^o \Bigl(e^{\la \kappa^L X_{\tau_1^{ (L)}} \cdot \ell} | 
\calF_0^0\Bigr)
= e^{\la L \kappa ^L }
 \oo{\EE}^o \Bigl(e^{\la \kappa ^L F} | \calF_0^0\Bigr) < Q(\la)<\infty\,, 
\; P-a.s.,
$$
with $\calF_0^0$ defined in (\ref{stronger}).
Further,
$\la_2\underset{\delta(\ell) \to 1}{\longrightarrow} 1$.
\end{lemma}
\medskip\noindent
\textbf{Proof of Lemma~\ref{lem-1}}

Set $W:=\min \{i:s_i \text{\ is L-successful}\}$. By definition, 
$s_W+L= \oo{S}_1$ and $X_{ \oo{S}_1}\cdot \ell = W+L$.
We first evaluate the
exponential moments of $W$.  
To every $i$, attach a random variable
$\chi_i = \won_{\{\oo\eps_{s_i+1} =e_1\}}$.  Note that the
$\chi_i$ are i.i.d., Bernoulli distributed with parameter $\kappa$ and
$$
W=\min \{j:\; \chi_j = \chi_{j+1} = \cdots = \chi_{j+L-1} = 1\}
\,.
$$
Consider  the inter-failure times  $\{\mu_i\}_{i \geq 1}$, i.e.
the sequence in $\NN^*$ defined by 
$$
\{j >0; \chi_j=0\}= \{ \mu_1, \mu_1+\mu_2, \ldots\}\;,
$$
which is i.i.d., geometrically distributed 
with failure probability $(1-\kappa)$,
and note that 
$W=\mu_1+ \ldots + \mu_i$ when 
$\mu_1 \leq L, \ldots, \mu_i \leq L,  \mu_{i+1} > L$.
 
With the notation $\oo{E}^o_\om= E_{Q } \oo{E}^o_{\om,{\beps}}$
as in the proof of Corollary \ref{cor-nonnest}, we have for all
$\omega$ such that the walk is $\oo{P}^o_\om$-a.s. transient in the
direction $\ell$,
\begin{align*}
\oo{E}^o_\om \exp \{ \la \kappa^L W \} &= 
\oo{E}^o_\om
 \sum_{i \geq 0} \exp \{ \la \kappa^L (\mu_1+ \ldots + \mu_i) \} \;
 {\bf 1}_{ \mu_1 < L, \ldots \mu_i < L, 
\mu_{i+1} \geq L }\\
 &= \sum_{i \geq 0}  \kappa^L \left( \oo{E}^o_\om
[{\bf 1}_{  \mu_1 < L}
\exp \{ \la \kappa^L \mu_1 \}] \right)^i\\
 &=
\frac{\kappa^L }{\left(1-\oo{E}^o_\om
[{\bf 1}_{  \mu_1 < L}
\exp \{ \la \kappa^L \mu_1 \}] \right)_+}\;,
\end{align*}
though
\begin{align*}
\oo{E}^o_\om[{\bf 1}_{  \mu_1 < L}\exp \{ \la \kappa^L \mu_1 \}]
= (1-\kappa) e^{ \la \kappa^L}  \frac{1- (\kappa e^{ \la \kappa^L})^L}
{1- (\kappa e^{ \la \kappa^L})} = 1 + (\la -1) \kappa^L + o(\kappa^L)\;.
\end{align*}
Hence, for all  $\la<1$, there exists a finite $L_1(\lambda)$ with
\begin{equation} \label{eq:4}
\sup_{L\geq L_1(\lambda)} {\rm ess \ sup}_{\om \in  {\rm  supp}(P)}
\oo{E}^o_\om 
\exp \{ \la \kappa^L W \} < \infty \;.
\end{equation}
Note that the estimate in \req{eq:4} is \textit{quenched\/}, i.e. for 
$P$-almost all environments.

Set 
$$
M_0= \max \Bigl\{
X_n \cdot \ell - X_{0} \cdot \ell\,,\;
0 \leq n <  D' \Bigr\} \in (0, \infty]\;,
$$
which is a.s. finite on the set $\{D'<\infty\}$ and infinite otherwise.
Next, let $\oo s$ denote an L-successful fresh time 
(i.e., $\oo s = s_k$ for some $k$),
and consider 
$$ 
M=M(\oo s):=  M_0 \circ \theta_{\oo s+L} = 
 \max \Bigl\{
X_n \cdot \ell - X_{\oo s+L} \cdot \ell\,,\;
\oo{s} + L \leq  n < \oo s + L + \theta_{\oo s+L} D' \Bigr\}\;.
$$  
Define 
$$
\oo\calF_{\oo s} = \sigma (\om_z : z\cdot \ell \le X_{\oo s}
 \cdot \ell) \vee \sigma (X_t, t \le \oo s)\;,
$$ 
recall the definition from (\ref{stronger})
${\cal F}_0^L=\sigma(\omega(y,\cdot): \ell\cdot y<-L)$ and the
notation $\bar \theta_x$ for space shift.
It is useful to  note that the 
annealed law of paths after fresh
points has the following  property,
\begin{align} \label{M=M_0loi}
\oo{\PP}^o (M > r | \oo{\calF}_{\oo s})=
\oo{\PP}^o (M_0 > r | {\cal F}_0^L) \circ \bar{\theta}_{X_{\oo{s} +L}} \;,\quad r >0\;,
\end{align}
which implies also that $\oo{\PP}^o (M > r, D'\circ \theta_{\oo{s} +L} < \infty
 | \oo{\calF}_{\oo s})=
\oo{\PP}^o (M_0 > r, D' < \infty | {\cal F}_0^L) \circ \bar{\theta}_{X_{\oo{s} +L}}$,
since $\{ M_0 >r, D'< \infty\}=\{  M_0 \in (r, \infty)\}$ almost surely.
Indeed, fix $n \geq L, \{x^*_i\}_{0 \leq i < n}$
a path for walk with $x^*_{i+1}=x^*_i+e_1$ for
$n-L \leq i \leq n$,
fix some $\om^* \in S^n$
and some measurable set $A \subset S^{\{y: y \cdot \ell \leq x_{n\!-\!L}^* \cdot
  \ell -L\}}$ such that $ \{\om (x^*_i)\}_{0 \leq i \leq n-L}= \om^* $
for all $\om \in A$, and $P(A)>0$. Then, by the 
Markov property, and for arbitrary
  $\om^{**} \in A$, 
\begin{align*}
\oo{\PP}^o (M > r &\; | \; {\oo s}=n-L, \{X_i\}_{0 \leq i \leq n}=
\{x^*_i\}_{0 \leq i \leq n} , \om \in A)\;=\\
&= \frac{ 
 E_{P\otimes Q}[
\oo{P}_{\om,\beps}^o (M > r, {\oo s}\! = \!n-L, \{X_i\}_{0 \leq i \leq n}\! = \!
\{x^*_i\}_{0 \leq i \leq n} ) {\bf 1}_{ \om \in A}]}
{\oo{\PP}^o ({\oo s}\! = \!n-L, \{X_i\}_{0 \leq i \leq n}\! = \!
\{x^*_i\}_{0 \leq i \leq n} , \om \in A)}\\
&= \frac{ 
 E_{P\otimes Q}
[\oo{P}_{\om, \theta_n \beps}^{x^*_n} (M_0 > r)
\oo{P}_{\om,\beps}^o ( {\oo s}\! = \!n-L, \{X_i\}_{0 \leq i \leq n}\! = \!
\{x^*_i\}_{0 \leq i \leq n}) {\bf 1}_{ \om \in A}]}
{\oo{\PP}^o ({\oo s}\! = \!n-L, \{X_i\}_{0 \leq i \leq n}\! = \!
\{x^*_i\}_{0 \leq i \leq n} , \om \in A)}\\
&= \frac{ 
 E_{P\otimes Q}
[\oo{P}_{\om, \theta_n \beps}^{x^*_n} (M_0 > r)
\oo{P}_{\om^{**},\beps}^o ( {\oo s}\! = \!n-L, \{X_i\}_{0 \leq i \leq n}\! = \!
\{x^*_i\}_{0 \leq i \leq n}) {\bf 1}_{ \om \in A}]}
{\oo{\PP}^o ({\oo s}\! = \!n-L, \{X_i\}_{0 \leq i \leq n}\! = \!
\{x^*_i\}_{0 \leq i \leq n} , \om \in A)}
\\
&= \frac{ 
 E_{ Q}[
E_P \left[
\oo{P}_{\om, \theta_n \beps}^{x^*_n} (M_0 > r)  {\bf 1}_{ \om \in A}
 \right]
\times
\oo{P}_{\om^{**},\beps}^o ( {\oo s}\! = \!n-L, \{X_i\}_{0 \leq i \leq n}\! = \!
\{x^*_i\}_{0 \leq i \leq n})]}
{ 
\oo{\PP}^o ({\oo s}\! = \!n-L, \{X_i\}_{0 \leq i \leq n}\! = \!
\{x^*_i\}_{0 \leq i \leq n} , \om \in A)}
\\
&= \frac{ 
 E_{P\otimes Q}\left[
\oo{P}_{\om, \theta_n \beps}^{x^*_n} (M_0 > r) {\bf 1}_{ \om \in A}
\right] \times
 E_{Q}\left[
\oo{P}_{\om^{**},\beps}^o ( {\oo s}\! = \!n-L, \{X_i\}_{0 \leq i \leq n}\! = \!
\{x^*_i\}_{0 \leq i \leq n})\right]}
{
E_{P\otimes Q}\left[
 {\bf 1}_{ \om \in A}
\right] \times
 E_{Q}\left[
\oo{P}_{\om^{**},\beps}^o ( {\oo s}\! = \!n-L, \{X_i\}_{0 \leq i \leq n}\! = \!
\{x^*_i\}_{0 \leq i \leq n})\right]}
\\
&=\frac{ 
 E_{P\otimes Q}
\oo{P}_{\om, \theta_n \beps}^{x^*_n} (M_0 > r, \om \in A )}{
P(\om \in A )}\;
=
\; \oo{\PP}^o (M_0 > r | \om \in \oo{\theta}_{x^*_n} A )
\end{align*}
performing the same computations on  the denominator, yielding (\ref{M=M_0loi}).

\begin{figure}[h]
\begin{center}
\includegraphics[scale=.6]{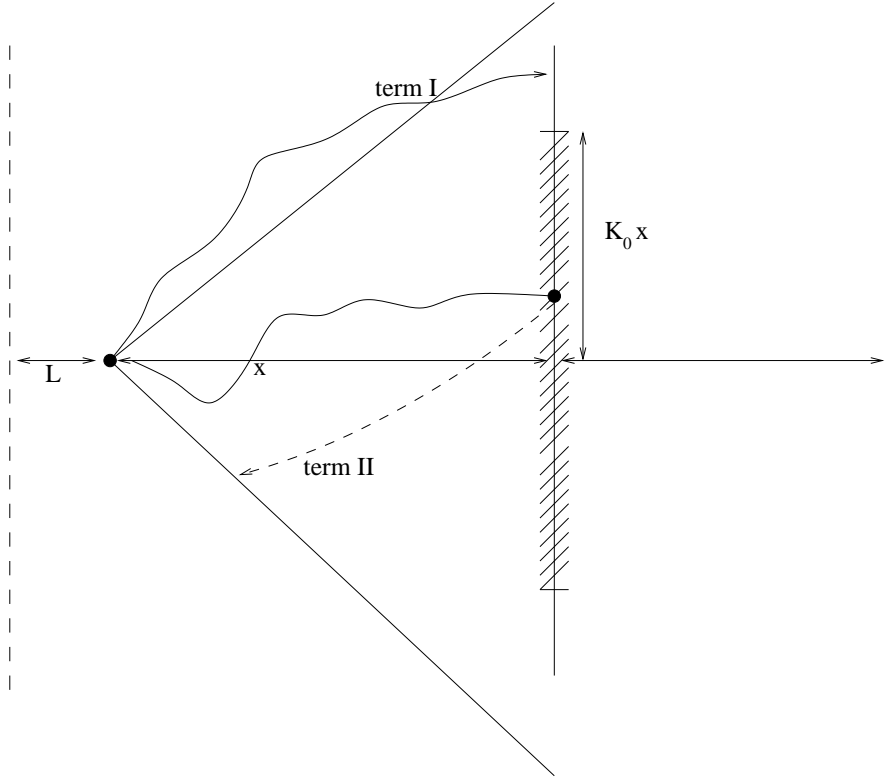}
\caption{Escape events}
\label{fig-2}
\end{center}
\end{figure}

We next derive tail estimates on $M$, considering only the case of $M_0$ in view
of (\ref{M=M_0loi}).
For $x \in \NN^*$, define
$$
U_x^{(L)} = \{ z\in \ZZ^d: -L \le z \cdot \ell \le x\}
$$
and set $\tau_{x,L} = \min \{n>0:\:X_{n} 
\not\in U_x^{(L)} \}$.
Finally, let $\Pi$ denote the orthogonal projection with respect to  $\ell=e_1$,
i.e.\ $\Pi (z) = z - z_1 e_1$, and $K_0$ a positive constant.
Then, for any $\zeta>0$, using the Markov property and a union bound,
\begin{align*}
\oo{\PP}^o (M_0 \in (x, \infty) \; | \; {\cal F}_0^L )
& \le \oo{\PP}^o \left( X_{\tau_{x,L}} \cdot
\ell \ge x,
\left| \Pi (X_{\tau_{x,L}} )
\right| > K_0 x \;| \; {\cal F}_0^L \right) \\
& + (2K_0 x)^{d-1} \max_{z: z \cdot \ell=x,
|\Pi(z)|\leq K_0x}
\;\oo{\PP}^z \Bigl( X_{\cdot} \,\mbox{\rm hits}
\,\, C(0,\ell,\zeta)^c\; | \; {\cal F}_0^L
\Bigr) 
=: \text{I+II}\;,
\end{align*}
as indicated in Figure \ref{fig-2}.

Note that the term I does not depend on $\zeta$.
We argue below that, 
by Kalikow's condition,
\begin{equation}
\label{osI}
\text{I} \le e^{-g(K_0,\delta(\ell)) x}
\end{equation}
where
$g(K_0,\delta(\ell)) \underset{K_0\to\infty}{\longrightarrow} \infty$
is monotone non decreasing in $\delta(\ell)$ and does not depend on $\zeta$.

Fix  now $K_0$ such that $g(K_0,\delta(\ell))>0$ and set $\zeta_0$ such that
\begin{equation}
\label{conelarge}
\{z: z\cdot\ell=x, |\Pi(z)|\leq 2K_0 x\}\subset C(0,\ell,\zeta_0)\}\,.
\end{equation}
We will also argue that
\begin{equation}
\label{osII}
\text{II} \le e^{- \beta_{K_0}' (\delta(\ell)) x}
\end{equation}
where $\beta_{K_0}'$ comes from Kalikow's condition and 
$\beta_{K_0}'(\delta) \underset{\delta \to 1}{\longrightarrow}\infty$.
Therefore,  
 we obtain for all $x$ large, and all $\zeta<\zeta_0$,
\begin{align} \label{tailM}
\sup_{\om} \oo{\PP}^o (M \in (x, \infty)  
 | \oo\calF_{\oo s}) \le  e^{-\beta (\delta
(\ell)
) x}
\end{align}
for some function $\beta=\beta(\ell)$ such that 
$ \beta \to \infty$ as $\delta \to 1$.
 
Indeed, to see (\ref{osI}) and thus control I, 
we use (\ref{sortie}) for the truncated strip $U_x^{(L)} \bigcap
\{ |\Pi(z)| \leq {\tilde L}\}$ for  large $\tilde L < \infty$, 
and apply the stopping theorem to the 
supermartingale $\exp\{ -3 \la_0 f(X_n)\}$ 
as in
Lemma  \ref{D'>0}-1)
(with the function 
$f(y)=y \cdot e_1 - \zeta' |y|, \zeta' =\delta(\ell)/2$),
 and to Kalikow's chain at 
the exit time from this strip. Letting then $\tilde L \to \infty$, 
one readily gets (\ref{osI}) 
with 
$g({K_0}) = 1.5 \la_0 (\delta(\ell) \sqrt{1+K_0^2}  -2)_+$. 
 To prove (\ref{osII}), we proceed similarly with truncated cones,
 using (\ref{sortie}) and, this time, the supermartingale 
$\exp\{ -3 \la_1 f(X_n)\}$
 as in 
 Lemma  \ref{D'>0}-4) (with the function 
$f(y)=y \cdot e_1 - \zeta_0 |y|$ and  assuming $ \zeta_0 \leq 
\delta(\ell)/3$ w.l.o.g.),
  to get
$$\oo{\PP}^z ( X_{\cdot} \,\mbox{\rm hits}\,\, C(0,\ell,\zeta)^c | 
 {\cal F}_0^L) \times \exp\{0\} \leq
\exp\{-3 \la_1x (1-  \zeta_0 \sqrt{K_0^2+1})\}$$
which yields (\ref{osII}) with $\beta_{K_0}'(\delta)= 
3 \la_1(\delta) (1-  \zeta_0 \sqrt{K_0^2+1})$, where the factor 
 $ 1-  \zeta_0 \sqrt{K_0^2+1} \geq 1-  (\sqrt{K_0^2+1}/\sqrt{4K_0^2+1}) >0$
from (\ref{conelarge}).

Next,  from Lemma \ref{D'>0}-4), it follows also that,
for any L-successful fresh time $\oo s$,
\begin{align} \label{Dlarge}
\oo{\PP}^o ( D' \circ \theta_{\oo{s}+L} =\infty | \oo{\calF}_{\oo s}) 
\geq^{^{(\ref{M=M_0loi})}}
q_1 := \inf_{\om} \oo{\PP}^o ( D' =\infty | \oo{\calF}_0^L) > 0
\end{align}
where
$q_1\underset{\delta\to 1}{\longrightarrow} 1$
uniformly in $L$ and $\zeta<\zeta_0$.

By definition, $X_{\oo{S}_1} \cdot e_1 = W + L$, and on $\{ \oo{S}_k < \infty\}$,
$ (X_{\oo{S}_{k+1}}-X_{\oo{S}_k}) \cdot e_1 = M_k + W_{k+1} + L$, with the notations
$M_k=M(\oo{S}_k)$, 
$W_{k+1}= W \circ \theta_{\oo{S}_k+M_k}$ if $M_k < \infty$, and
$W_{k+1}= \infty$ on $\{M_k = \infty \}= \{D' \circ \theta_{\oo{S}_k}
= \infty\}$.
\begin{figure}[h]
\begin{center}
\includegraphics{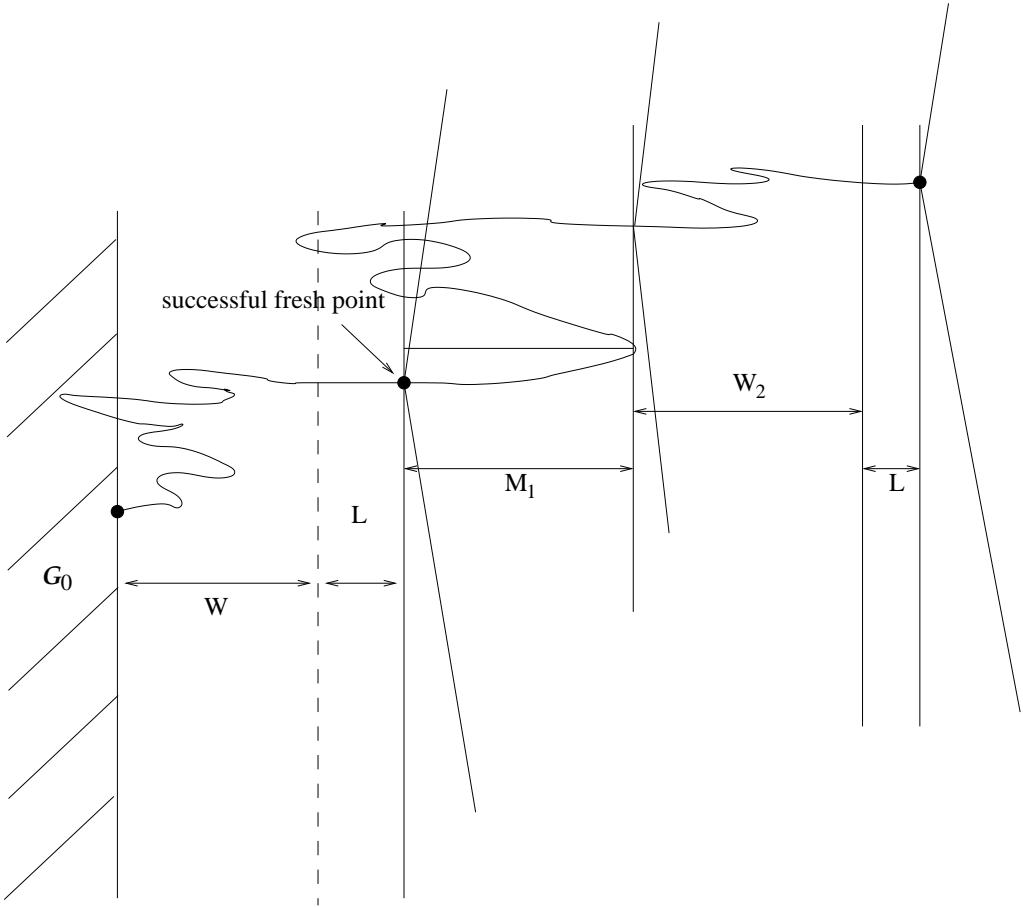}
\caption{Successful fresh 
times leading to regeneration.
}
\label{fig-3}
\end{center}
\end{figure}
Put also $W_1=W$, write 
$$
F+L 
= \sum_{k \geq 0} (M_k {\bf 1}_{k \neq 0}+ W_{k+1} + L) 
{\bf 1}_{k < K}\;,
$$
(see Figure \ref{fig-3}) and, with $\la >0$, 
\begin{eqnarray*}
&&  \!\!  \!\! \! \! 
\oo{\EE}^o(  
 e^{\la \kappa^L (F+L)} \;|\;  \calF_0^0)  \leq 
\oo{\EE}^o ( \sum_{k \geq 0}
 e^{\la \kappa^L \sum_{l=0}^k (M_l {\bf 1}_{l \neq 0}
 + W_{l+1} + L) }
{\bf 1}_{M_1, \ldots M_k< \infty}
\;|\; \calF_0^0)\\
= &&\!\!\!  \!\!\!
\sum_{k \geq 0} 
\oo{\EE}^o \left( e^{\la \kappa^L \sum_{l=0}^{k-1}(M_l {\bf 1}_{l \neq 0}
+ W_{l+1} + L) } 
{\bf 1}_{M_1, \ldots M_{k-1}< \infty}\!\! \times \!\!
\oo{\EE}^o\left[ e^{\la \kappa^L (M_k {\bf 1}_{k \neq 0}
+ W_{k+1} + L) } {\bf 1}_{M_k< \infty}
| \oo{ \calF}_{\oo{S}_k}\vee {\calF}_0^0\right] \;
|\; \calF_0^0 \right)\\
\leq  && \!\!\! \!\! \!
\sum_{k \geq 0} 
\oo{\EE}^o \left( e^{\la \kappa^L \sum_{l=0}^{k-1}(M_l {\bf 1}_{l \neq 0}
 + W_{l+1} + L) } 
{\bf 1}_{M_1, \ldots M_{k-1}< \infty}| \calF_0^0 \right)  \!\!\times \!\!
\oo{\EE}^o[e^{\la \kappa^L (W+L)}] \sup_{\om }
\oo{\EE}^o [ e^{\la \kappa^L M_0 
} {\bf 1}_{D'< \infty}
|  {\calF}_0^0] \\
\leq & &\!\!\! \!\! \!
\sum_{k \geq 0} \left( 
\oo{\EE}^o[e^{\la \kappa^L (W+L)}]\right)^{k+1} \left(\sup_{\om }
\oo{\EE}^o [ e^{\la \kappa^L M_0 } {\bf 1}_{D'< \infty}
|  {\calF}_0^0]\right)^{k} \\
= &&\!\!\!\!\!\!
\frac{e^{\la L \kappa^L } \oo{\EE}^o[e^{\la \kappa^L W)}]}{\Big(1-
  e^{\la L \kappa^L }
\oo{\EE}^o[e^{\la \kappa^L W}]
\sup_{\om } \oo{\EE}^o [ e^{\la \kappa^L M_0 } {\bf 1}_{D'< \infty}
|  {\calF}_0^0]\Big)_+} =: Q(\la )\;,
\end{eqnarray*}
where the third line comes from independence and from (\ref{M=M_0loi}), and
 the fourth one 
from a recursion. From (\ref{eq:4}), $\sup_{ L> L_1(\lambda)} \{
e^{\la L \kappa^L } \oo{\EE}^o[e^{\la \kappa^L W}]\}< \infty$ for $\la <1$, while 
by Schwarz's inequality,
$$
\sup_{\om } \oo{\EE}^o [ e^{\la \kappa^L M_0 } {\bf 1}_{D'< \infty}
|  {\calF}_0^0] \leq \left[(1-q_1) \sup_{\om }
 \oo{\EE}^o [ e^{2\la \kappa^L M_0 }{\bf 1}_{M_0<\infty}
| {\calF}_0^0] 
\right]^{1/2}\;.
$$
From (\ref{tailM}), we can choose
 $L$ large enough so that the supremum in the right-hand 
side is finite, and therefore there exists some $\la_2 >0$ 
such that $Q( \la )< \infty$ for
$\la < \la_2$. This completes the proof of the first claim in the Lemma.
To obtain the second claim, we  make $1-q_1$ arbitrary small by taking
$\delta( \ell )$ close to 1, thus keeping  $Q$ finite for $\la$ arbitrarily
 close to 1.

This completes the proof of Lemma \ref{lem-1}.
\qed

We now complete the proof of Theorem \ref{theo-A5}, by deriving
 probability estimates for tails of $\tau_1^{(L)}$.
Fix $n>0$, $T=3 \kappa^{-L}  \ln n$, and a large constant ${K_1}>0$.  Set
\begin{align*}
B & = \{ z\in \ZZ^d: \; -L \le  z\cdot\ell \le T \}\,,\\
\tilde B_{K_1} & = \{ z\in B:\; \:|\:\Pi(z)\:|\: \le {K_1}T\}, \; 
\partial_+ \tilde B_{K_1}=
\{z\in \partial \tilde B_{K_1}:\; z\cdot \ell \ge T\}
\,.
\end{align*}
(Recall that $\Pi$ is the projection on the hyperplane orthogonal to
$\ell$.) With $T_{B^c}$ 
[resp. $T_{{\tilde B}_{K_1}^c}$,] the 
exit time from ${B}$ [resp. ${\tilde {B}_{K_1}}$],
we decompose the set $\{\tau_1^{(L)} > \kappa^{-L} n\} \cap \{D'=\infty\}$ 
according to
$X_{\cdot}$ exiting  $\tilde B_{K_1}$ after time $\kappa^{-L} n$, and 
then decompose 
 the latter case according to $X_{\cdot}$ exiting  $\tilde B_{K_1}$ 
from $\partial_+\tilde B_{K_1}$ or not:
\begin{eqnarray}
\label{eq:10}
 \oo{\PP}^o \Bigl(\tau_1^{(L)} 
 > \kappa^{-L} n\;,\;  D'=\infty \:|\: \calF_{0}^L\Bigr)
& \leq &  \oo{\PP}^o \Bigl( {T_{{\tilde B}_{K_1}^c}}> \kappa^{-L} n\;,\; D'=\infty
\:|\: \calF_{0}^L\Bigr)
 + \oo{\PP}^o \Bigl( X_{\tau_1^{ (L)}}
\cdot \ell > T\:|\: 
\calF_{0}^L \Bigr)  + \notag\\
&& + \oo{\PP}^o \Bigl(  X_{
T_{{\tilde B}_{K_1}^c}
} \notin \partial_+\tilde B_{K_1}\;, \;
D'=\infty  \:|\: \calF_{0}^L \Bigr) \notag\\
&
\le &
 \oo{\PP}^o
\Bigl(  {T_{{\tilde B}_{K_1}^c}}> \kappa^{-L} n\;,\; D'=\infty
\:|\: \calF_{0}^L\Bigr) + e^{-3 \la \ln n} 
+ 0\;,
\end{eqnarray}
where $\la < \la_2$ from Lemma \ref{lem-1}, and where the last term is
made equal to zero by fixing  $K_1 > \zeta^{-1}(1- \zeta^2)^{1/2}$, in which
case all path contained in $C(0, \ell, \zeta)$ have to exit $\tilde
B_{K_1}$ from $\partial_+\tilde B_{K_1}$. 

In $\tilde B_{K_1}$, consider strips of width $\Delta = \delta_1\ln n$, with
$-\delta_1\ln\kappa<1$ and $T/\Delta$ integer,  
$$
B_i = \{z\in \tilde B_{K_1} :\: (i-1) \Delta < z \cdot \ell < (i+1)
\Delta\}\;,
\quad i =0,1, \cdots,
T/\Delta
$$
consider the truncated hyperplanes (slices)
$A_i = \{ z\in \tilde B_{K_1}:\: z\cdot \ell = [i \Delta] 
\}$, and define the random variables ($T_{B^c_i}=$ exit
time of $B_i$)
$$
Y_i:=\sup_{z\in A_i} P_\om^z ( X_{T_{B^c_i}}\cdot \ell < (i-1) \Delta)\;,
$$
i.e., the smallest quenched probability starting from the middle slice
$A_i$ to exit the strip $B_i$ from the left.

By Kalikow's condition, for $z \in A_i$,
$$
 E\Bigl( P_\om^z(X_{T_{B^c_i}}\cdot \ell < (i-1) \Delta) \;\vert\; \om(x,
 \cdot), x \notin B_i
\Bigr) \le e^{-c \Delta}
$$
with $c=c(\delta)\underset{\delta\to1}{\longrightarrow} \infty$
(using the supermartingale $\exp\{-2 \la_1 f(X_n)\}$ from Lemma
\ref{D'>0}, see proof of (\ref{tailM})). 
Hence,
$$
 P\Bigl( Y_i > e^{-c\Delta/2} \;
\vert \;\om(x, \cdot), x \notin B_i\Bigr)
 \le (2{K_1}T)^{d-1}\; e^{-c\Delta/2} 
$$
In particular, the set 
$$
{\calA}:=\{ \exists i\leq T/\Delta :\: Y_i > e^{-c\Delta/2}\}\;
\text{ is such that }\;
P\Bigl( {\calA} \; |\;
\calF_0^L \Bigr) \le
\frac{T}{\Delta } (2{K_1}T)^{d-1}\; e^{-c\Delta/2} 
\,.
$$
Let now ${\cal B}$ denote the event that the walk, sampled at hitting times
of neighboring slices, successively visits $A_1,
A_2, \dots A_{T/\Delta }$ without backtracking to the neighboring
slice on the right. 
Then, for $\om \in \calA^c$,
$$
P_\om^0 \left({\cal B}^c  \right)
\le 1- \Bigl(1-e^{-c\Delta/2}\Bigr)^{\ffrac{T}{\Delta} } 
\le \: \frac{T}{\Delta} e^{- c\Delta/2} \,,
$$
by convexity ($ T/\Delta>1$).
On the other hand, the times spent inside each block are stochastically
dominated above by $2(\Delta +1)  G_i$, where
$G_i$ are independent  random variables, geometrically distributed with parameter
$\kappa^{2(\Delta +1)}$; indeed, by ellipticity, the walk starting from
any point in $B_i$ has probability larger than $\kappa^{2(\Delta +1)}$
to exit  $B_i$ by traveling only with steps to the right, with at most 
$2(\Delta +1)$ steps.
Now, since
$$
 \oo{\PP}^o
\Bigl(  {T_{{\tilde B}_{K_1}^c}}\!>\! \kappa^{-L} n\;,\; D'\!=\!\infty
\:|\: \calF_{0}^L\Bigr) \leq 
 \oo{\PP}^o \Bigl({\cal A} \:|\: \calF_{0}^L \Bigr) + 
 \oo{\PP}^o \Bigl( {\cal A}^c \bigcap  {\cal B}^c \:|\: \calF_{0}^L \Bigr) + 
E \Bigl( P_\om^o(  {T_{{\tilde B}_{K_1}^c}}\!>\! \kappa^{-L} n\;,\; {\cal B})
{\bf 1}_{ {\cal A}^c} \:|\: \calF_{0}^L \Bigr)\;,
$$
we have
\begin{align*}
& \oo{\PP}^o \Bigl( \{\tau_1^{(L)} > \kappa^{-L} n\} \cap \{ D'=\infty\} \:|\:
\calF_{0}^L\Bigr)  \\
& \qquad \le \left[n^{-3 \la}
 + 
\frac{T}{\Delta } (2{K_1}T)^{d-1}\; e^{-c\Delta/2} 
 + \frac{T}{\Delta} e^{-c\Delta/2} 
+ P \left( \sum_{i \leq T/\Delta}
 G_i > \frac{\kappa^{-L} n}{2(\Delta +1)} \right)
\right] \wedge 1 
\,.
\end{align*}
But, on the last event, at least one of the $G_i$'s is larger than 
${ \Delta \kappa^{-L} n}/{[2(\Delta +1)T]} $, so
\begin{align*}
P \left( \sum_{i \leq T/\Delta}
 G_i > \frac{\kappa^{-L} n}{2(\Delta +1)} \right)
&  \le 
\frac{T}{\Delta} P\left(G_1 > \frac{n}{12 \ln n } \right) \\
&  = \frac{T}{\Delta} \Bigl(1-\kappa^{2(\Delta +1)}
\Bigr)^{\frac{n}{12 \ln n}}  \leq 
\frac{3 \kappa^{-L}}{\delta_1} 
\exp\{- \frac{n}{12 n^{- 2 \delta_1 \ln \kappa} \ln n }\}\,,
\end{align*}
and finally
\begin{align*}
& \oo{\PP}^o \Bigl( \{\tau_1^{(L)} > \kappa^{-L} n\} \cap \{ D'=\infty\} \:|\:
\calF_{0}^L \Bigr)  \\
& \qquad \qquad \le \left[ n^{-3 \la}
 + \frac{3 \kappa^{-L} }{ \delta_1}
\left( n^{-\delta_1 c/2} 
[(6{K_1} \kappa^{-L} \ln n)^{d-1}+1 ] + e^{- \frac{n}{12 n^{-
      2 \delta_1 
\ln \kappa} \ln n }}\right) \right]\wedge 1 \\
& \qquad \qquad
\le \begin{cases}
1, & \mbox{\rm for} \; \;K_2 \kappa^{-dL} (\ln n)^{d-1}
> n^{+\delta_1 c/4}\;, \\ \\
2 n^{-[(\delta_1 c/4)\vee (2 \la)]}, & \mbox{\rm else}\,,
\end{cases}
\end{align*}
for some constant $K_2$,  all $n$ large enough and all $\la < \la_2$.
Hence, with $\alpha=2$,
$$
M(L)= 
\EE^o\Bigl( (\tau_1^{(L)} \kappa^{L})^2 1_{\{D=\infty\}}| {\cal F}_{0}^L\Bigr)
\le K_3  \kappa^{-8dL/c\delta_1}\,,
$$
for some constant $K_3$ independent of $L$, as soon as $c=c(\delta)$ is large 
enough, which happens 
if $\kappa$ is kept fixed and $\delta\to 1$. 
Thus, as soon as $\phi'(\cdot)$ decreases exponentially 
one may
find a 
$\delta(\ell)$
close enough to $1$ such that $M(L)\phi'(L)
\to_{L\to\infty} 0\,.$ For such $\delta(\ell)$, we thus conclude that
$({\cal A}5)$ is satisfied with $\alpha=2$.
\qed

We are now ready to describe the class of examples
satisfying the assumptions of Theorem \ref{theo-A5}.
Let $(\sigma(x), x \in \ZZ^d)$ 
be a $d$-dimensional nearest 
neighbor Ising model with $\beta \geq 0, h>0$ and $\pi$ the
corresponding probability measure, i.e. 
$$
\pi( \sigma(x)= \pm 1 | \sigma(y), y \neq x) 
= \exp  \pm \{\beta {\tt S} +h\}
/(\exp   \{\beta {\tt S} +h\}+ \exp - \{\beta {\tt S} +h\})\;,
$$
with
${\tt S}=\sum_{e: |e|=1}
\sigma(x+e)$.
Fix now two probability vectors 
$\om^{\pm}=(\om^{\pm}(e);  e\in \ZZ^d,
|e|=1)$, and assume that 
\begin{equation} \label{ell-ex}
\om^{\pm}(e)>0 \;,\;\; |e|=1\;,
\end{equation}
\begin{equation} \label{drift-ex}
d^+:= 
\left(\sum_{|e|=1}
\om^+(e)e \right) \cdot e_1 >0\;,\quad 
-d^-:= 
\left(\sum_{|e|=1}
\om^-(e)e \right) \cdot e_1 <0\;.
\end{equation}
Consider the random environment given by 
\begin{equation} \label{def-ex}
\om(x,x+e)= \om^{\pm}(e) \quad {\rm according \ to  \ } \; \sigma(x)=\pm 1\;.
\end{equation}
Note that the RWRE is nestling in the case where the local drifts
points in opposite directions, e.g.
$\sum_{|e|=1} \om^+(e)e \in (0,1) \cdot e_1$ and 
$\sum_{|e|=1} \om^-(e)e \in (-1,0) \cdot e_1$. 
\begin{theorem}
\label{theo-exgibbs} 
For all choice of $\om^{\pm}$ 
with (\ref{ell-ex}) and (\ref{drift-ex}),
there exist 
a finite number $h_0$
and a positive function $\beta_0(h)$ with 
$\lim_{h \to + \infty } \beta_0(h)=
\infty$ such that for $h>h_0$ and $\beta<\beta_0(h)$,
$$\lim_{n\to\infty}
\frac1n\sum_{i=1}^n X_n=v\,,\quad \oo{\PP}^o-a.s..$$
for some deterministic vector $v \;(v \cdot e_1 >0)$.
\end{theorem}

Note that, since the above function $\beta_0$ is unbounded, 
the result applies to
arbitrary low temperature - but with large external field.

\noindent
{\bf Proof of Theorem \ref{theo-exgibbs}}
We start by giving  a sufficient condition for Kalikow's condition 
$({\cal A}3)$.
\begin{lemma}
\label{lemma-CSKala} 
If for a 
deterministic $\delta >0$ we have $P$-a.s.
\begin{align} \label{CSKal}
\inf_{ f: \{\pm e_i\}_1^d \to  (0,1]} 
\left[  E\left( \frac{1}
{\sum_{|e|=1} f(e) \om (x,x\!+\!e)} 
 | {\calF}_{\{x\}^c} \right)^{-1}
 E\left( \frac{d(x, \om )\cdot \ell}
{\sum_{|e|=1} f(e) \om (x,x\!+\!e)}
| {\calF}_{\{x\}^c} \right)
\right] \geq \delta \,,
\end{align}
then $({\calA} 3)$ holds with $\ell$
and $\delta(\ell)=\delta$. 
\end{lemma}

The proof is similar  to \cite{kalikow},
p.759-760, replacing 
$\mu$ therein with our Gibbs measure 
$ \pi(.|\sigma(y), y \neq x)$. We will omit it here,
and we will use Lemma \ref{lemma-CSKala} with $\ell=e_1$.

\begin{lemma}
\label{lemma-CSKal}
For condition (\ref{CSKal}) to hold with $\delta >0$ and $\ell=e_1$
in the example (\ref{def-ex}),  it is sufficient that
$$
  \delta < d^+\;,\quad 2h -4 \beta d \geq \ln 
\frac{(d^-+\delta)}{d^+-\delta} +\ln 
\max_{|e|=1} \frac
{\om^+ (e)}{\om^- (e)}
\;.
$$
\end{lemma}

\noindent
{\bf Proof of Lemma \ref{lemma-CSKal}}
To simplify notations, we set ${\tilde \pi}( \pm | {\tt S})=
\pi( \sigma(x)= \pm 1 | \sigma(y), y \neq x)$ on the set $\{\sum_{|e|=1}
\sigma(x+e)=s\}$, and ${\cal S}=\{-2d, -2d+2, \ldots, 2d \}$ 
the set of possible 
values for  $s$. For $\delta >0$,  
\begin{eqnarray}
(\ref{CSKal}) \quad & \Longleftrightarrow & \forall f, s, \quad
\frac
{
 {\tilde \pi}( + | s)
  \frac{d^+}{\sum_{|e|=1} f(e) \om^+ (e)}
 +{\tilde \pi}( - | s)
  \frac{-d^-}{\sum_{|e|=1} f(e) \om^- (e)} 
}{
 {\tilde \pi}( + | s)
  \frac{1}{\sum_{|e|=1} f(e) \om^+ (e)}
 +{\tilde \pi}( - | s)
  \frac{1}{\sum_{|e|=1} f(e) \om^- (e)} 
} \geq \delta \nonumber \\
& \Longleftrightarrow & \forall f, s,\quad
\frac
{ {\tilde \pi}( + | s)}{{\tilde \pi}( - | s)}
 \geq 
\frac
{\sum_{|e|=1} f(e) \om^+ (e)}{\sum_{|e|=1} f(e) \om^- (e)} 
\times 
\frac
{d^-+\delta}{d^+ -\delta}\;,
\label{d5}
\end{eqnarray}
since $\delta < d^+$. Note that for the Ising measure $\pi$,
$$
\inf_{s \in {\cal S}} \frac
{ {\tilde \pi}( + | s)}{{\tilde \pi}( - | s)}
= \inf_{s \in {\cal S}} \exp(2 h + 2\beta s)
\geq \exp(2 h - 4 \beta d)\;,
$$
while on the other hand,
$$
\sup_{ f: \{\pm e_i\}_1^d  \to (0,1]} \frac
{\sum_{|e|=1} f(e) \om^+ (e)}{\sum_{|e|=1} f(e) \om^- (e)} =
\max_{|e|=1} \frac
{\om^+ (e)}{\om^- (e)}\;.
$$
Lemma \ref{lemma-CSKal}
is proved. \qed

Let $c=c(\beta, h)$ be Dobrushin's contraction coefficient (for
example, definition (2.7) in
\cite{follmer}). 
If $c=c(h, \beta)<1$ then $\pi$ is weak-mixing, 
with a constant $\gamma$ depending only on $c$ - as can be checked 
from (2.8) in \cite{follmer}. (See also \cite[Theorem 3, Section 2.2.1.3]{doukhan}.)
According to 
Proposition \ref{d51},
 $P$ satisfies $({\cal A}1)$ with $\phi(L)\leq e^{-\gamma' L}$, where
 $\gamma'=\gamma /\sqrt{2}$ depends only on $c$. From
Theorem \ref{theo-A5} and Lemma \ref{lemma-CSKal}, Condition  
$({\cal A}5)$ is implied by
$$
c(h, \beta)<1\;,\;\; 
  \delta < d^+\;,\;\; 
2h -4 \beta d \geq \ln 
\frac{(d^-+\delta)}{d^+-\delta} +\ln 
\max_{|e|=1} \frac
{\om^+ (e)}{\om^- (e)}
\;,\;\; 
\delta > \delta_0( \om^-(e_1), \gamma'(c(h, \beta)),d)\;.
$$
Note that $c(h, \beta)\leq c_0<1$ 
contains a region $0 \leq \beta < \beta_0'(h)$ with $\beta_0'(h) \to \infty$
as $h \to \infty$. The proof of Theorem \ref{theo-exgibbs}
is complete. \qed

\section{Concluding remarks}
\label{sec:conclusions}
\setcounter{equation}{0}
\begin{enumerate}
\item While working on the first draft of this work, we learnt of 
\cite{komorowski}, where the authors consider law of large numbers
for $L$-dependent non-nestling environments. The approach 
of \cite{komorowski} is quite different from ours, as it relies on
constructing an invariant measure, absolutely continuous
with respect to the law $P$, which makes the regeneration
sequence $\{\bar \tau_i^{(0)}, \bar X_i^{(0)}\}$ (with $L=0$ and
$\zeta =0$ as  introduced in \cite{sznitmanzerner})
a stationary sequence. 
These results are covered by our approach.

In another preprint \cite{komorowski2}, the same authors obtain 
the  law of large numbers for $L$-dependent
  environments with Kalikow's condition,
under the crucial assumption that with positive 
probability the walker jumps at a 
distance larger than $L$. An appropriate modification of the regeneration
times introduced in \cite{sznitmanzerner} leads in that case
to a renewal structure as in the 
independent environment.

We also mention that \cite{firas} has announced results related to ours, 
obtained by the method of the environment viewed from the point 
of view of the article.

\item It is reasonable to expect that under $({\cal A}1,2,3)$,
the integrability 
condition $({\cal A}5)$, with $\alpha>2$, implies 
that the CLT for $X_n/\sqrt{n}$ holds true. However,
the proof of such a statement does present some challenges. We hope to
return to this question in future work.
\item It is worthwhile to note that in the case of i.i.d. 
environments, under Kalikow's condition
it holds that $({\cal A}5)$ with any 
$\alpha>1$ is satisfied as soon as $d\geq 2$, by the results of 
\cite{Sz1}. It is not clear whether, under reasonable mixing 
conditions which are not of the $L$-dependent type, the law of large numbers 
holds for the whole range of Kalikow's condition, or more generally
what should the analogue of Sznitman's T'-condition (see \cite{Sz2})
be in the mixing setup.
\end{enumerate}

\end{document}